\documentclass[final,leqno]{siamltex}

\usepackage[leqno]{amsmath}
\usepackage{graphicx}
\usepackage{graphics}
\usepackage{epstopdf}
\usepackage{subfigure}
\usepackage{threeparttable}

\usepackage{bm}
\usepackage{amssymb}
\usepackage{leftidx}
\usepackage{empheq}
\usepackage{color}
\usepackage{makecell}
\allowdisplaybreaks
\renewcommand{\theequation}{\arabic{section}.\arabic{equation}}
\numberwithin{equation}{section}
\newtheorem{remark}{Remark}[section]
\title{Step-by-step solving schemes based on scalar auxiliary variable and invariant energy quadratization approaches for gradient flows.
        \thanks{
We would like to acknowledge the assistance of volunteers in putting together this example manuscript and supplement. This work is supported by the Postdoctoral Science Foundation of China under grant numbers BX20190187 and 2019M650152, by National Natural Science Foundation of China (Grant Nos: 11901489, 11971276).}}

      \author{Zhengguang Liu
             \thanks{School of Mathematics and Statistics, Shandong Normal University, Jinan, China. Email: liuzhgsdu@yahoo.com}.
                                       \and
             Xiaoli Li\textsuperscript{*}
             \thanks{Corresponding author: Fujian Provincial Key Laboratory on Mathematical Modeling and High Performance Scientific Computing and School of Mathematical Sciences, Xiamen University, Xiamen, Fujian, 361005, China. Email: xiaolisdu@163.com}. }
\begin{document}

\maketitle

\begin{abstract}
In this paper, we propose several novel numerical techniques to deal with nonlinear terms in gradient flows. These step-by-step solving schemes, termed 3S-SAV and 3S-IEQ schemes, are based on recently popular scalar auxiliary variable (SAV) and invariant energy quadratization (IEQ) approaches. In these constructed numerical methods, the phase function $\phi$ and auxiliary variable can be calculated step-by-step. Compared with the traditional SAV/IEQ approaches, there are many advantages for the novel 3S-SAV/3S-IEQ schemes. Firstly, we do not need the restriction of the bounded from below of the nonlinear free energy potential/density function. Secondly, the auxiliary variable combined with nonlinear function can be treated totally explicitly in the 3S-SAV/3S-IEQ approaches. Specially, for solving the discrete scheme based on IEQ approach, the linear system usually involves variable coefficients which change at each time step. However, the discrete scheme based on 3S-IEQ approach leads to linear equation with constant coefficients. Two comparative studies of traditional SAV/IEQ and 3S-SAV/3S-IEQ approaches are considered to show the accuracy and efficiency. Finally, we present various 2D numerical simulations to demonstrate the stability and accuracy.
\end{abstract}

\begin{keywords}
Step-by-step solving scheme, scalar auxiliary variable, invariant energy quadratization, gradient flows, numerical simulations.
\end{keywords}

    \begin{AMS}
         65M12; 35K20; 35K35; 35K55; 65Z05.
    \end{AMS}

\pagestyle{myheadings}
\thispagestyle{plain}
\markboth{ZHENGGUANG LIU AND XIAOLI LI} {3S-SAV/3S-IEQ APPROACHES FOR GRADIENT FLOWS}
  \section{Introduction}
The gradient flows are very important equations in physics, material science and mathematics \cite{ambati2015review,guo2015thermodynamically,liu2019efficient,marth2016margination,miehe2010phase,shen2015efficient,wheeler1992phase,wheeler1993computation}. They have been widely used in many fields such as alloy casting, new material preparation, image processing, finance and so on. As we all know, the phase field models are significant and important parts of gradient flows. They were proposed as mathematical models to describe an isothermal material phase transition of a binary alloy. Nowadays, phase field methodology can approximate various moving interface problems. Specially, many physical phenomena, such as the formation process of snowflakes, the dendrite structure formed by water freezing, and the cellular or dendrite structure formed in the welding process, etc can be simulated by gradient flows. It is helpful to understand the nature and the formation mechanism of various materials and the preparation of new materials. It is of great practical significance to develop new technologies.

Mathematically, the gradient flow is derived from the functional variation of free energy. In general, the free energy $E(\phi)$ contains the sum of an integral phase of a nonlinear functional and a quadratic term:
\begin{equation}\label{intro-e1}
E(\phi)=\frac12(\phi,\mathcal{L}\phi)+E_1(\phi)=\frac12(\phi,\mathcal{L}\phi)+\int_\Omega F(\phi)d\textbf{x},
\end{equation}
where $\mathcal{L}$ is a symmetric non-negative linear operator, and $E_1(\phi)$ is nonlinear free energy. $F(x)$ is the energy density function. The gradient flow from the energetic variation of the above energy functional $E(\phi)$ in \eqref{intro-e1} can be obtained as follows:
\begin{equation}\label{intro-e2}
  \left\{
   \begin{array}{rll}
\displaystyle\frac{\partial \phi}{\partial t}&=&\mathcal{G}\mu,\\
\mu&=&\displaystyle\mathcal{L}\phi+F'(\phi),
   \end{array}
   \right.
\end{equation}
where $\mu=\frac{\delta E}{\delta \phi}$ is the chemical potential. $\mathcal{G}$ is a non-positive operator. $\mathcal{G}=-I$ for the Allen-Cahn type system and $\mathcal{G}=\Delta$ for the Cahn-Hilliard type system for the Ginzburg-Landau double-well type potential $F=\frac{1}{4}(\phi^2-1)^2$.

It is not difficult to find that the gradient flow system satisfies the following energy dissipation law:
\begin{equation*}
\frac{d}{dt}E=(\frac{\delta E}{\delta \phi},\frac{\partial\phi}{\partial t})=(\mathcal{G}\mu,\mu)\leq0,
\end{equation*}
which is a very important property for gradient flows in physics and mathematics. In numerical computation, energy dissipation law becomes a very important indicator for numerical schemes. Up to now, many scholars considered a series of efficient and popular time discretized approaches to construct energy stable schemes for gradient flows such as convex splitting approach \cite{eyre1998unconditionally,shen2012second,shin2016first}, linear stabilized approach \cite{shen2010numerical,yang2017numerical}, exponential time differencing (ETD) approach \cite{du2019maximum,WangEfficient}, invariant energy quadratization (IEQ) approach \cite{chen2019efficient,chen2019fast,yang2016linear}, scalar auxiliary variable (SAV) approach \cite{xiaoli2019energy,shen2018scalar,ShenA} and so on. There are many different advantages for these methods. Specifically, the convex splitting method leads to a convex minimization problem at each time step and the scheme is unconditionally energy stable and uniquely solvable. But it still needs to solve a nonlinear system, and it is difficult to construct higher-order scheme. The linear stabilized method can effectively solve a linear system of the gradient flows, but the stabilized term leads to additional error, which makes it difficult to construct the higher-order scheme. The IEQ method is inspired by Lagrange multiplier method but makes a big leap. By introducing a auxiliary variable, X. Yang and his collaborators \cite{YangNumerical,YangEfficient} successfully avoided the difficulty of discretization of the nonlinear term. The IEQ approach has been proven to keep many advantages such as linear, easy to obtain second order scheme and unconditional energy stability. This method has been successfully applied to the numerical simulation for many complex phase field models. The SAV method was proposed by J. Shen and his collaborators \cite{shen2018scalar,ShenA} which is another very popular and efficient approach. It is worth mentioning that the SAV method keeps all the advantages of the IEQ approach. Furthermore, it weakens the assumptions of the bounded below restriction of nonlinear free energy potential.

In recent years, both IEQ and SAV approaches have become very efficient and powerful ways to construct energy stable schemes for gradient flows. However, both of them still have some obvious shortcomings. We list some of the issues here:

$(i)$ the assumption conditions that the nonlinear free energy $E_1(\phi)$/density function $F(\phi)$ in SAV/IEQ approaches must be bounded from below are required to keep the square root reasonable. It means that there exists a constant $C$ to satisfy $E_1(\phi)+C>0$/$F(\phi)+C>0$. However, it is observed that the minimum values of $E_1(\phi)$ and $F(\phi)$ are not easy to obtain before calculation. Furthermore, the positive constant $C$ seems to have an influence on the accuracy of the simulation results \cite{lin2019energy,lin2019numerical}.

$(ii)$ one can not guarantee the positive property of the introducing auxiliary variables $r=\sqrt{E_1(\phi)+C}$ and $q=\sqrt{F(\phi)+C}$ in SAV/SEQ approaches, to say noting of the numerical $r^n$ and $q^n$ in numerical schemes. If the positive property of $r$ can not be guaranteed in SAV approach, the term $\frac{r}{\sqrt{E_1(\phi)+C}}$ might not be 1.

$(iii)$ the totally explicit schemes of the auxiliary variable combined with the nonlinear term are not energy stable for SAV approach. Then, an inner product has to be calculated before obtaining $\phi$. Such procedure requires us to solve equations to obtain the values of inner products for the multiple SAV approach.

$(iv)$ for solving the discrete scheme based on IEQ approach, although one only needs to solve a linear system at each time step, the linear system usually involves variable coefficients which change at each time step.

In order to enhance the applicability of IEQ and SAV approaches and improve the above mentioned issues, we propose some novel numerical techniques to deal with nonlinear terms. The novel auxiliary variables $r$ and $q$ which do not include square root are introduced in the proposed 3S-SAV and 3S-IEQ approaches. The favorable properties include:

$(i)$ the proposed 3S-SAV and 3S-IEQ approaches weaken the assumption conditions. We do not need a unknown constant $C$ to keep $E_1(\phi)+C>0$ and $F(\phi)+C>0$. We just need $E_1(\phi)+C\neq0$ and $F(\phi)+C\neq0$. It can be proved that such constant $C$ is very easy to obtain before calculation.

$(ii)$ the auxiliary variables $r$ and $q$ in the proposed 3S-SAV and 3S-IEQ approaches will not include square root which are more reasonable for the equivalence of the models.

$(iii)$ the auxiliary variable combined with the nonlinear term can be treated totally explicitly in the 3S-SAV/3S-IEQ approaches. In calculation, the phase function $\phi$ and auxiliary variable can be calculated step-by step. Specially, the phase function $\phi$ be solved directly rather than computing the inner product.

$(iv)$ the discrete scheme based on 3S-IEQ approach leads to linear equations with constant coefficients so it is remarkably easy to implement.

The paper is organized as follows. In Sect.2, to show a comparative study for our modified approaches, we give a brief review of IEQ and SAV approaches.  Then, we introduce the 3S-SAV approach and consider first and second order time discrete schemes for gradient flows in Sect.3. All discrete schemes are proved the unconditional energy stability. In Sect.4, a step-by-step solver based on IEQ approach is considered. Finally, in Sect.5, various 2D numerical simulations are demonstrated to verify the accuracy and efficiency of our proposed schemes.
\section{A brief review of IEQ and SAV approaches}
In order to show and give a comparative study for our modified approaches, we provide below a brief review of the popular IEQ and SAV approaches to construct energy stable schemes for gradient flows. Before giving a semi-discrete formulation, we let $N>0$ be a positive integer and set
\begin{equation*}
\Delta t=T/N,\quad t^n=n\Delta t,\quad \text{for}\quad n\leq N.
\end{equation*}
\subsection{Invariant energy quadratization (IEQ) approach}
The core idea of IEQ approach is to transform the nonlinear potential in gradient flows into a simple quadratic form. Then the nonlinear term $F'(\phi)$ can be treated explicitly. What's more, the derivative of the quadratic polynomial is linear, one only needs to solve the linear equations with constant coefficients at each time step. In detail, assume that the energy density function $F(\phi)$ is bounded from below which means that there is a constant $C>0$ to satisfy $F(\phi)+C>0$. Then, introduce an auxiliary function as follows
\begin{equation*}
q(x,t;\phi)=\sqrt{F(\phi)+C}.
\end{equation*}
The nonlinear term $F'(\phi)$ can be transformed as follows:
\begin{equation*}
F'(\phi)=\frac{q}{q}F'(\phi)=\frac{q}{\sqrt{F(\phi)+C}}F'(\phi).
\end{equation*}

Using above transformation and noting that $q_t=\frac{1}{2\sqrt{F(\phi)+C}}F'(\phi)\phi_t$, we can rewrite the gradient flows \eqref{intro-e2} as the following equivalent system:
\begin{equation}\label{ieq-e1}
  \left\{
   \begin{array}{rll}
\displaystyle\frac{\partial \phi}{\partial t}&=&\mathcal{G}\mu,\\
\mu&=&\displaystyle\mathcal{L}\phi+\frac{q}{\sqrt{F(\phi)+C}}F'(\phi),\\
q_t&=&\displaystyle\frac{F'(\phi)}{2\sqrt{F(\phi)+C}}\phi_t.
   \end{array}
   \right.
\end{equation}

Taking the inner products of the above equations with $\mu$, $\phi_t$ and $2q$, respectively, it is easy to obtain that the above equivalent system \eqref{ieq-e1} satisfies a modified energy dissipation law:
\begin{equation*}
\frac{d}{dt}\left[\frac12(\phi,\mathcal{L}\phi)+\int_{\Omega}q^2d\textbf{x}\right]=(\mathcal{G}\mu,\mu)\leq0.
\end{equation*}

Based on the principle of dealing with nonlinear terms explicitly and linear terms implicitly, we give a second-order scheme based on the Crank-Nicolson method easily:
\begin{equation}\label{ieq-e2}
  \left\{
   \begin{array}{rll}
\displaystyle\frac{\phi^{n+1}-\phi^{n}}{\Delta t}&=&\mathcal{G}\mu^{n+1/2},\\
\mu^{n+1/2}&=&\displaystyle\mathcal{L}\left(\frac{\phi^{n+1}+\phi^n}{2}\right)+\frac{q^{n+1}+q^n}{2\sqrt{F(\tilde{\phi}^{n+1/2})+C}}F'(\tilde{\phi}^{n+1/2}),\\
\displaystyle\frac{q^{n+1}-q^n}{\Delta t}&=&\displaystyle\frac{1}{2\sqrt{F(\tilde{\phi}^{n+1/2})+C}}F'(\tilde{\phi}^{n+1/2})\frac{\phi^{n+1}-\phi^{n}}{\Delta t},
   \end{array}
   \right.
\end{equation}
where $\tilde{\phi}^{n+\frac{1}{2}}$ is any second order explicit approximation for $\phi(t^{n+\frac{1}{2}})$, which can be flexible according to the problem.

It is not difficult to prove that above scheme is unconditionally energy stable in the sense that
\begin{equation*}
\aligned
\left[\frac12(\mathcal{L}\phi^{n+1},\phi^{n+1})+\int_{\Omega}(q^{n+1})^2d\textbf{x}\right]-\left[\frac12(\mathcal{L}\phi^{n},\phi^{n})+\int_{\Omega}(q^n)^2d\textbf{x}\right]\leq\Delta t(\mathcal{G}\mu^{n+1/2},\mu^{n+1/2})\leq0.
\endaligned
\end{equation*}
\subsection{Scalar auxiliary variable (SAV) approach}
The SAV approach is an enhanced version of the IEQ approach. That is because that the SAV method keeps all the advantages of the IEQ approach. Furthermore, it weakens the assumptions of the bounded below restriction of nonlinear free energy potential which makes it to be a new important way to simulate the gradient flows. The key of the SAV approach is to transform the nonlinear potential $E_1(\phi)$ into a simple scalar quadratic form. In particulat, assuming that $E_1(\phi)$ is bounded from below which means that there exists a constant $C$ to make $E_1(\phi)+C>0$. Define a scalar auxiliary variable
\begin{equation*}
r(t)=\sqrt{E_1(\phi)+C}=\sqrt{\int_\Omega F(\phi)d\textbf{x}+C}>0.
\end{equation*}
Then, the nonlinear functional $F'(\phi)$ can be transformed into the following equivalent formulation:
\begin{equation*}
F'(\phi)=\frac{r}{r}F'(\phi)=\frac{r}{\sqrt{E_1(\phi)+C}}F'(\phi).
\end{equation*}
Thus, an equivalent system of gradient flow \eqref{intro-e2} with scalar auxiliary variable can be rewritten as follows
\begin{equation}\label{sav-e1}
  \left\{
   \begin{array}{rll}
\displaystyle\frac{\partial \phi}{\partial t}&=&\mathcal{G}\mu,\\
\mu&=&\displaystyle\mathcal{L}\phi+\frac{r}{\sqrt{E_1(\phi)+C}}F'(\phi),\\
r_t&=&\displaystyle\frac{1}{2\sqrt{E_1(\phi)+C}}\int_{\Omega}F'(\phi)\phi_td\textbf{x}.
   \end{array}
   \right.
\end{equation}

Taking the inner products of the above equations with $\mu$, $\phi_t$ and $2r$, respectively, we obtain that the above equivalent system satisfies a modified energy dissipation law:
\begin{equation*}
\frac{d}{dt}\left[\frac12(\phi,\mathcal{L}\phi)+r^2\right]=(\mathcal{G}\mu,\mu)\leq0.
\end{equation*}

The above equivalent system \eqref{sav-e1} is very easy to construct linear, second order and unconditional energy stable scheme. For example, a second-order semi-discrete scheme based on the Crank-Nicolson method reads as follows
\begin{equation}\label{sav-e2}
  \left\{
   \begin{array}{rll}
\displaystyle\frac{\phi^{n+1}-\phi^{n}}{\Delta t}&=&\mathcal{G}\mu^{n+1/2},\\
\mu^{n+1/2}&=&\displaystyle\mathcal{L}\left(\frac{\phi^{n+1}+\phi^n}{2}\right)+\frac{r^{n+1}+r^n}{2\sqrt{E_1(\tilde{\phi}^{n+1/2})+C}}F'(\tilde{\phi}^{n+1/2}),\\
\displaystyle\frac{r^{n+1}-r^n}{\Delta t}&=&\displaystyle\frac{1}{2\sqrt{E_1(\tilde{\phi}^{n+1/2})+C}}\int_{\Omega}F'(\tilde{\phi}^{n+1/2})\frac{\phi^{n+1}-\phi^{n}}{\Delta t}d\textbf{x}.
   \end{array}
   \right.
\end{equation}

It is not difficult to obtain the following discrete dissipation law:
\begin{equation*}
\aligned
\left[\frac12(\mathcal{L}\phi^{n+1},\phi^{n+1})+|r^{n+1}|^2\right]-\left[\frac12(\mathcal{L}\phi^{n},\phi^{n})+|r^{n}|^2\right]\leq\Delta t(\mathcal{G}\mu^{n+1/2},\mu^{n+1/2})\leq0.
\endaligned
\end{equation*}
No restriction on the time step is required.
\section{3S-SAV approach for gradient flow}
In this section, we will give our novel 3S-SAV approach and describe its advantages in calculation in detail. Firstly, we give an explicit scheme for gradient flows based on the exactly same auxiliary variable with the traditional SAV approach. This explicit scheme can be proved to be unconditional energy stable. However, such scheme is not suitable for computation, since the term in square root can not be proved to keep positive. By introducing a novel scalar auxiliary variable without square root to replace the former one, we obtain the 3S-SAV approach for gradient flow.
\subsection{3S-SAV approach based on traditional auxiliary variable}
In this section, we will consider a 3S-SAV approach to construct an explicit scheme with unconditionally energy stability. In particular, using the exactly same auxiliary variable $r(t)$ with the introduced SAV approach in above section:
\begin{equation}\label{3s-sav}
r(t)=\sqrt{E_1(\phi)+C}, \quad F'(\phi)=\frac{r}{\sqrt{E_1(\phi)+C}}F'(\phi).
\end{equation}

To simplify the notations, we define a new variable $\chi^{r,\phi}$ as follows:
\begin{equation*}
\chi^{r,\phi}=\frac{r}{\sqrt{E_1(\phi)+C}}F'(\phi).
\end{equation*}

The system \eqref{sav-e1} based on SAV approach can be rewritten as the following equivalent system:
\begin{equation}\label{3s-sav-e1}
  \left\{
   \begin{array}{rll}
\displaystyle\frac{\partial \phi}{\partial t}&=&\mathcal{G}\mu,\\
\mu&=&\displaystyle\mathcal{L}\phi+\chi^{r,\phi},\\
\displaystyle2r\frac{dr}{dt}&=&\displaystyle(\chi^{r,\phi},\frac{\partial \phi}{\partial t}).
   \end{array}
   \right.
\end{equation}

Taking the inner products of the first two equations with $\mu$ and $\phi_t$ in \eqref{3s-sav-e1} respectively, then combining them with the third equation in \eqref{3s-sav-e1} and noting that $2rr_t=\frac{dr^2}{dt}$, we obtain the following energy dissipation law:
\begin{equation*}
\frac{d}{dt}\left[\frac12(\mathcal{L}\phi,\phi)+r^2\right]=(\mathcal{G}\mu,\mu)\leq0.
\end{equation*}

It is easy to find that the above energy inequality is totally equal to the energy dissipation law based on the SAV approach.

Next, we will find that the system \eqref{3s-sav-e1} is very easy to construct linear and unconditionally energy stable semi-discrete schemes.

\subsubsection{The first-order scheme}
A first order scheme for solving the system \eqref{3s-sav-e1} can be readily derived by the backward Euler¡¯s method. The first-order scheme can be written as follows:
\begin{equation}\label{3s-sav-first-e1}
  \left\{
   \begin{array}{rll}
\displaystyle\frac{\phi^{n+1}-\phi^{n}}{\Delta t}&=&\mathcal{G}\mu^{n+1},\\
\mu^{n+1}&=&\displaystyle\mathcal{L}\phi^{n+1}+\chi^{r^n,\phi^n},\\
\displaystyle(r^{n+1}+r^n)\frac{r^{n+1}-r^n}{\Delta t}&=&\displaystyle\left(\chi^{r^n,\phi^n},\frac{\phi^{n+1}-\phi^{n}}{\Delta t}\right),
   \end{array}
   \right.
\end{equation}
Multiplying the first two equations in \eqref{3s-sav-first-e1} with $\mu^{n+1}$ and $(\phi^{n+1}-\phi^{n})/\Delta t$, and combining  them with the third equation in \eqref{3s-sav-first-e1}, we immediately obtain the discrete energy law:
\begin{equation}\label{3ssav-first-e2}
\aligned
\frac{1}{\Delta t}\left[E_{1st}^{n+1}-E^{n}_{1st}\right]\leq(\mathcal{G}\mu^{n+1},\mu^{n+1})\leq0,
\endaligned
\end{equation}
where the modified discrete version of the energy is defined by
\begin{equation*}
\aligned
E_{1st}^{n}=\frac12(\phi^n,\mathcal{L}\phi^{n})+(r^n)^2.
\endaligned
\end{equation*}

The proposed first-order 3S-SAV scheme \eqref{3s-sav-first-e1} is much easier to calculate than traditional SAV. For the SAV scheme which can be seen in \cite{ShenA}, to obtain $\phi^{n+1}$, we have to compute an inner product $(b^n,\phi^{n+1})$ previously where $b^n=F'(\phi^n)/\sqrt{E_1(\phi^n)}$. However, for 3S-SAV scheme, we only need to compute $\phi^{n+1}$ and $r^{n+1}$ step-by-step. Because $\phi^{n+1}$ can be calculated directly by the first two equations in \eqref{3s-sav-first-e1}, then $r^{n+1}$ can be very easy to obtain by computing $\left(b^{r^n,\phi^n},\phi^{n+1}-\phi^{n}\right)$. Thus compared with the SAV algorithm, the 3S-SAV algorithm greatly simplifies the calculation which is conducive for rapid simulation.

Particularly, we rewrite the first two equation in \eqref{3s-sav-first-e1} as the following matrix formulation:
\begin{equation}\label{3s-sav-first-e3}
\aligned
(I-\Delta t\mathcal{G}\mathcal{L})\phi^{n+1}=\phi^n+\Delta t\mathcal{G}\chi^{r^n,\phi^n}.
\endaligned
\end{equation}
Multiplying \eqref{3s-sav-first-e3} with $(I-\Delta t\mathcal{G}\mathcal{L})^{-1}$, we can obtain $\phi^{n+1}$ directly:
\begin{equation}\label{3s-sav-first-e4}
\aligned
\phi^{n+1}=(I-\Delta t\mathcal{G}\mathcal{L})^{-1}\phi^n+\Delta t(I-\Delta t\mathcal{G}\mathcal{L})^{-1}\mathcal{G}\chi^{r^n,\phi^n}.
\endaligned
\end{equation}
Substitute equation \eqref{3s-sav-first-e4} into the third equation in \eqref{3s-sav-first-e3}, we can compute $r^{n+1}$:
\begin{equation}\label{3s-sav-first-e5}
r^{n+1}=\displaystyle\sqrt{(r^n)^2+\left(\chi^{r^n,\phi^n},\phi^{n+1}-\phi^{n}\right)}.
\end{equation}

\begin{remark}\label{3s-sav-re1}
From above analysis, we find that the computation of $\phi^n$ and $r^n$ can be solved step by step.  However,the scheme  \eqref{3s-sav-first-e3}-\eqref{3s-sav-first-e5} may be blow up due to the fact that one can not guarantee positive property for  $(r^n)^2+\left(b^{r^n,\phi^n},\phi^{n+1}-\phi^{n}\right)$ in square root in \eqref{3s-sav-first-e5}.
\end{remark}

Next, to overcome this difficulty, we rewrite the nonlinear term $F'(\phi)$ in \eqref{3s-sav} as follows:
\begin{equation}\label{3s-sav-first-e6}
\aligned
F'(\phi)=\displaystyle\frac{r^2}{r^2}F'(\phi)=\displaystyle\frac{r^2}{E_1(\phi)+C}F'(\phi).
\endaligned
\end{equation}
Then, we redefine the variable $\chi^{r,\phi}$ as follows:
\begin{equation*}
\chi^{r,\phi}=\frac{r^2}{E_1(\phi)+C}F'(\phi).
\end{equation*}
Then, the system \eqref{3s-sav-e1} can be transformed as follows:
\begin{equation}\label{3s-sav-e2}
  \left\{
   \begin{array}{rll}
\displaystyle\frac{\partial \phi}{\partial t}&=&\mathcal{G}\mu,\\
\mu&=&\displaystyle\mathcal{L}\phi+\chi^{r,\phi},\\
\displaystyle\frac{dr^2}{dt}&=&\displaystyle(\chi^{r,\phi},\frac{\partial \phi}{\partial t}).
   \end{array}
   \right.
\end{equation}

Then, the first-order scheme can be written as follows:
\begin{equation}\label{3s-sav-e3}
  \left\{
   \begin{array}{rll}
\displaystyle\frac{\phi^{n+1}-\phi^{n}}{\Delta t}&=&\mathcal{G}\mu^{n+1},\\
\mu^{n+1}&=&\displaystyle\mathcal{L}\phi^{n+1}+\chi^{r^n,\phi^n},\\
\displaystyle\frac{(r^{n+1})^2-(r^n)^2}{\Delta t}&=&\displaystyle\left(\chi^{r^n,\phi^n},\frac{\phi^{n+1}-\phi^{n}}{\Delta t}\right),\\
\chi^{r^n,\phi^n}&=&\displaystyle\frac{(r^n)^2}{E_1(\phi^n)+C}F'(\phi^n).
   \end{array}
   \right.
\end{equation}

It is easy to obtain that the discrete scheme \eqref{3s-sav-e3} has the following discrete energy law:
\begin{equation}\label{3s-sav-e4}
\aligned
\frac{1}{\Delta t}\left[\left(\frac12(\phi^n,\mathcal{L}\phi^{n+1})+(r^{n+1})^2\right)-\left(\frac12(\phi^n,\mathcal{L}\phi^{n})+(r^n)^2\right)\right]\leq(\mathcal{G}\mu^{n+1},\mu^{n+1})\leq0.
\endaligned
\end{equation}

Noting that $\chi^{r^n,\phi^n}=\displaystyle\frac{(r^n)^2}{E_1(\phi^n)+C}F'(\phi^n).$ It means that we only need to compute $(r^{n+1})^2$ rather than to solve the value of $r^{n+1}$. Particularly, $(r^{n+1})^2$ can be obtained as follows:
\begin{equation}\label{3s-sav-e5}
(r^{n+1})^2=\displaystyle(r^n)^2+\left(\chi^{r^n,\phi^n},\phi^{n+1}-\phi^{n}\right).
\end{equation}
\subsection{The first-order scheme based on 3S-SAV approach with novel SAV}
In this subsection, we try to give a novel scalar auxiliary variable to re-derive the system \eqref{3s-sav-e2} and the discrete numerical scheme \eqref{3s-sav-e3}. Introduce a scalar auxiliary variable $\eta(t)$:
\begin{equation}\label{3s-sav-e6}
\eta(t)=E_1(\phi)+C=\int_\Omega F(\phi)d\textbf{x}+C\neq0.
\end{equation}
Then, the nonlinear functional $F'(\phi)$ can be transformed into the following equivalent formulation:
\begin{equation*}
F'(\phi)=\frac{\eta}{\eta}F'(\phi)=\frac{\eta}{E_1(\phi)+C}F'(\phi).
\end{equation*}

Taking the derivative of $\eta(t)$ with respect to $t$ and using above transformation of $F'(\phi)$, we obtain
\begin{equation}\label{3s-sav-e7}
\frac{d\eta(t)}{dt}=\frac{dE_1(\phi)}{dt}=\int_\Omega F'(\phi)\phi_td\textbf{x}=\int_\Omega \frac{\eta}{E_1(\phi)+C}F'(\phi)\phi_td\textbf{x}.
\end{equation}

To simplify the notations, we define a variable $\chi^{\eta,\phi}$ as follows:
\begin{equation*}
\chi^{\eta,\phi}=\frac{\eta}{E_1(\phi)+C}F'(\phi).
\end{equation*}

Then, the gradient flow \eqref{intro-e2} can be transformed into the following equivalent formulation:
\begin{equation}\label{3s-sav-e8}
  \left\{
   \begin{array}{rll}
\displaystyle\frac{\partial \phi}{\partial t}&=&\mathcal{G}\mu,\\
\mu&=&\displaystyle\mathcal{L}\phi+\chi^{\eta,\phi},\\
\displaystyle\frac{d\eta}{dt}&=&\displaystyle(\chi^{\eta,\phi},\frac{\partial \phi}{\partial t}).
   \end{array}
   \right.
\end{equation}
We can easily obtain a modified energy dissipation law by taking the inner products of the above first two equation with $\mu$ and $\phi_t$. It reads
\begin{equation*}
\frac{d}{dt}\left(\frac12(\mathcal{L}\phi,\phi)+\eta\right)=(\mathcal{G}\mu,\mu)\leq0.
\end{equation*}

Similarly, the first-order scheme based on 3S-SAV approach can be written as follows:
\begin{equation}\label{3s-sav-e9}
  \left\{
   \begin{array}{rll}
\displaystyle\frac{\phi^{n+1}-\phi^{n}}{\Delta t}&=&\mathcal{G}\mu^{n+1},\\
\mu^{n+1}&=&\displaystyle\mathcal{L}\phi^{n+1}+\chi^{\eta^n,\phi^n},\\
\displaystyle\frac{\eta^{n+1}-\eta^n}{\Delta t}&=&\displaystyle\left(\chi^{\eta^n,\phi^n},\frac{\phi^{n+1}-\phi^{n}}{\Delta t}\right),\\
\chi^{\eta^n,\phi^n}&=&\displaystyle\frac{\eta^n}{E_1(\phi^n)+C}F'(\phi^n).
   \end{array}
   \right.
\end{equation}
It is easy to obtain that the discrete scheme \eqref{3s-sav-e3} has the following discrete energy law which no restriction on the time step is required.:
\begin{equation}\label{3s-sav-e10}
\aligned
\frac{1}{\Delta t}\left[\left(\frac12(\phi^n,\mathcal{L}\phi^{n+1})+\eta^{n+1}\right)-\left(\frac12(\phi^n,\mathcal{L}\phi^{n})+\eta^n\right)\right]\leq(\mathcal{G}\mu^{n+1},\mu^{n+1})\leq0.
\endaligned
\end{equation}

In 3S-SAV scheme \eqref{3s-sav-e9}, the computations of $\phi$ and the auxiliary variable $\eta$ are totally decoupled. We do not calculate the inner product before obtaining $\phi$. We can compute $\phi^n$ and $\eta^n$ by a step-by-step solver. Secondly, we introduce a new auxiliary variable $\eta$ with no square root. What we should point out is that one should choose a proper constant $C$ carefully to satisfy $E_1(\phi^n)+C\neq0$.

\begin{remark}\label{3s-sav-re2}
In 3S-SAV scheme \eqref{3s-sav-e9}, we need to choose a proper constant $C$ to satisfy $E_1(\phi^n)+C\neq0$. For gradient flows, the dissipative energy law means $\frac{d}{dt}E(\phi)\leq0$. Then, an obvious property will hold as follows
\begin{equation*}
E(\phi(\textbf{x},0))\geq E(\phi(\textbf{x},t)), \quad \forall \textbf{x}\in\Omega,t\geq0.
\end{equation*}
Considering the definition of the energy in \eqref{intro-e1} and noting that $\mathcal{L}$ is a symmetric non-negative linear operator, it is not difficult to obtain the following inequality
\begin{equation*}
E(\phi(\textbf{x},0))-E_1(\phi(\textbf{x},t))=E(\phi(\textbf{x},0))-E(\phi(\textbf{x},t))+(\phi,\mathcal{L}\phi)\geq(\phi,\mathcal{L}\phi)\geq0, \quad \forall \textbf{x}\in\Omega,t\geq0.
\end{equation*}
It means that we can choose $C=-E(\phi(\textbf{x},0))-\delta$ where $\delta$ is an arbitrary positive constant.
\end{remark}

\begin{remark}\label{3s-sav-re3}
The 3S-SAV scheme \eqref{3s-sav-e9} can be solved efficiently by using a step-by-step solver as follows:\\
$~~~~~~(i)$ Compute $\eta^0=E_1(\phi_0)+C$ by using the initial condition $\phi(x,0)=\phi_0$;\\
$~~~~~~(ii)$ Compute $\chi^{\eta^n,\phi^n}$ from $\displaystyle\chi^{\eta^n,\phi^n}=\frac{\eta^n}{\sqrt{E_1(\phi^n)+C}}F'(\phi^n)$ for $n\geq0$;\\
$~~~~~~(iii)$ Compute $\phi^{n+1}$ from $\phi^{n+1}=(I-\Delta t\mathcal{G}\mathcal{L})^{-1}\phi^n+\delta t(I-\Delta t\mathcal{G}\mathcal{L})^{-1}\mathcal{G}\chi^{\eta^n,\phi^n}$;\\
$~~~~~~(iv)$ Compute $\eta^{n+1}$ from $\eta^{n+1}=\displaystyle\eta^n+\left(\chi^{\eta^n,\phi^n},\phi^{n+1}-\phi^{n}\right)$;\\
$~~~~~~(v)$ Let $n=n+1$, and go back to step $(ii)$.
\end{remark}
\subsection{The second-order scheme}
A linear, second-order, sequentially solved and unconditionally stable 3S-SAV scheme is also very easy to construct.  In detail, a semi-implicit 3S-SAV scheme based on the second order Crank-Nicolson formula (CN) for \eqref{3s-sav-e8} reads as: for $n\geq1$,
\begin{equation}\label{3s-sav-second-e1}
  \left\{
   \begin{array}{rll}
\displaystyle\frac{\phi^{n+1}-\phi^{n}}{\Delta t}&=&\mathcal{G}\mu^{n+\frac12},\\
\mu^{n+\frac12}&=&\displaystyle\mathcal{L}\frac{\phi^{n+1}+\phi^{n}}{2}+\chi^{\widetilde{\eta}^{n+\frac12},\widetilde{\phi}^{n+\frac12}},\\
\displaystyle\frac{\eta^{n+1}-\eta^n}{\Delta t}&=&\displaystyle\left(\chi^{\widetilde{\eta}^{n+\frac12},\widetilde{\phi}^{n+\frac12}},\frac{\phi^{n+1}-\phi^{n}}{\Delta t}\right),\\
\chi^{\widetilde{\eta}^{n+\frac12},\widetilde{\phi}^{n+\frac12}}&=&\displaystyle\frac{\widetilde{\eta}^{n+\frac12}}{E_1(\widetilde{\phi}^{n+\frac12})-E(\phi_0)-\delta}F'(\widetilde{\phi}^{n+\frac12}),
   \end{array}
   \right.
\end{equation}
where $\tilde{\phi}^{n+\frac{1}{2}}$ is any explicit $O(\Delta t^2)$ approximation for $\phi(t^{n+\frac{1}{2}})$, and $\tilde{\chi}^{n+\frac{1}{2}}$ is any explicit $O(\Delta t^2)$ approximation for $\chi(t^{n+\frac{1}{2}})$, which can be flexible according to the problem. Here, we choose
\begin{equation}\label{3s-sav-second-e2}
\aligned
&\tilde{\phi}^{n+\frac{1}{2}}=\frac32\phi^n-\frac12\phi^{n-1}, \quad n\geq1,\\
&\tilde{\chi}^{n+\frac{1}{2}}=\frac32\chi^n-\frac12\chi^{n-1}, \quad n\geq1,
\endaligned
\end{equation}
and for $n=0$, we compute $\widetilde{\phi}^{\frac{1}{2}}$ and $\widetilde{\chi}^{\frac{1}{2}}$ as follows:
\begin{equation}\label{3s-sav-second-e3}
\aligned
&\displaystyle\frac{\widetilde{\phi}^{\frac{1}{2}}-\phi^0}{(\Delta t)/2}=\mathcal{G}\left[\mathcal{L}\widetilde{\phi}^{\frac{1}{2}}+F^{'}(\phi^0)\right],\\
&\displaystyle\widetilde{\chi}^{\frac{1}{2}}=\int_\Omega F(\widetilde{\phi}^{\frac{1}{2}})d\textbf{x},
\endaligned
\end{equation}
which have local truncation errors of $O(\Delta t^2)$.

Similarly, $\phi^{n+1}$ and $\eta^{n+1}$ can be solved step-by-step:
\begin{equation}\label{3s-sav-second-e5}
\aligned
&\phi^{n+1}=(I-\frac12\Delta t\mathcal{G}\mathcal{L})^{-1}\phi^n+\frac12(I-\frac12\Delta t\mathcal{G}\mathcal{L})^{-1}\mathcal{G}\mathcal{L}\phi^n+\Delta t(I-\Delta t\mathcal{G}\mathcal{L})^{-1}\mathcal{G}\chi^{\widetilde{\eta}^{n+\frac12},\widetilde{\phi}^{n+\frac12}},\\
&\eta^{n+1}=\displaystyle\eta^n+\left(\chi^{\widetilde{\eta}^{n+\frac12},\widetilde{\phi}^{n+\frac12}},\phi^{n+1}-\phi^{n}\right).
\endaligned
\end{equation}

Multiplying the first two equations in \eqref{3s-sav-second-e1} with $\mu^{n+\frac12}$ and $(\phi^{n+1}-\phi^{n})/\Delta t$, and combining  them with the third equation in \eqref{3s-sav-second-e1}, we derive the following theorem immediately:
\begin{theorem}\label{3s-sav-th1}
The scheme \eqref{3s-sav-second-e1} for the equivalent gradient flow system \eqref{3s-sav-e8} is second-order accurate, unconditionally energy stable in the sense that
\begin{equation*}
\aligned
\frac{1}{\Delta t}\left[E_{3S-SAV/CN}^{n+1}-E^{n}_{3S-SAV/CN}\right]\leq(\mathcal{G}\mu^{n+\frac12},\mu^{n+\frac12})\leq0.
\endaligned
\end{equation*}
where the modified discrete version of the energy is defined by
\begin{equation*}
\aligned
E_{3S-SAV/CN}^{n}=\frac12(\phi^n,\mathcal{L}\phi^{n})+\eta^n.
\endaligned
\end{equation*}
\end{theorem}
\begin{remark}\label{esav-re1}
The second-order 3S-SAV scheme \eqref{3s-sav-second-e1} based on Crank-Nicolson can be implemented step-by-step as follows: (i) Compute the initial values of $\phi^0$ and $\eta^0$; (ii) Calculate $\tilde{\phi}^{\frac{1}{2}}$ and $\tilde{\eta}^{\frac{1}{2}}$ by using \eqref{3s-sav-second-e3}; (iii) Calculate $\phi^1$ by using  the first equation in \eqref{3s-sav-second-e5}; (iv) Calculate $\eta^1$ by using the second equation in \eqref{3s-sav-second-e5}; (v) Calculate $\tilde{\phi}^{n+\frac{1}{2}}$ and $\tilde{\eta}^{n+\frac{1}{2}}$ from \eqref{3s-sav-second-e2} for $n\geq1$; (vi) Calculate $\phi^{n+1}$ and $\eta^{n+1}$ from \eqref{3s-sav-second-e5}.
\end{remark}
\section{3S-IEQ approach for gradient flows}
Similar to 3S-SAV approach, it is easy to construct 3S-IEQ approach for gradient flows. In this section, by introducing a novel Lagrange multiplier $q(x,t;\phi)$, we give a step-by step solver based on IEQ approach. Specifically, assuming that there exists an constant $C$ to make $F(\phi)+C\neq0$. Then, introduce an auxiliary variable
\begin{equation}\label{3s-ieq-e1}
q(x,t;\phi)=F(\phi)+C.
\end{equation}

We rewrite the nonlinear term $F'(\phi)$ as the following formulation:
\begin{equation}\label{3s-ieq-e2}
F'(\phi)=\frac{q}{q}F'(\phi)=\frac{q}{F(\phi)+C}F'(\phi).
\end{equation}

Then, the gradient flow system \eqref{intro-e2} can be transformed into the following:
\begin{equation}\label{3s-ieq-e3}
  \left\{
   \begin{array}{rll}
\displaystyle\frac{\partial \phi}{\partial t}&=&\mathcal{G}\mu,\\
\mu&=&\displaystyle\mathcal{L}\phi+\chi^{q,\phi},\\
\displaystyle\frac{dq}{dt}&=&\displaystyle\chi^{q,\phi}\phi_t,\\
\chi^{q,\phi}&=&\displaystyle\frac{q}{F(\phi)+C}F'(\phi).
   \end{array}
   \right.
\end{equation}

Taking the inner products of the above equations with $\mu$, $\phi_t$, and $1$ respectively, we obtain the following  modified energy dissipation law:
\begin{equation*}
\frac{d}{dt}\left[\frac12(\phi,\mathcal{L}\phi)+\int_{\Omega}qd\textbf{x}\right]=(\mathcal{G}\mu,\mu)\leq0.
\end{equation*}

Both first and second-order discrete schemes with unconditionally energy stability can be obtained immediately. For example, the second-order 3S-IEQ scheme can be written $n\geq2$:
\begin{equation}\label{3s-ieq-e5}
  \left\{
   \begin{array}{rll}
\displaystyle\frac{\phi^{n+1}-\phi^{n}}{\Delta t}&=&\mathcal{G}\mu^{n+\frac12},\\
\mu^{n+\frac12}&=&\displaystyle\mathcal{L}\frac{\phi^{n+1}+\phi^{n}}{2}+\chi^{q^{n},q^{n-1},\phi^{n},\phi^{n-1}},\\
\displaystyle\frac{q^{n+1}-q^n}{\Delta t}&=&\displaystyle\chi^{q^{n},q^{n-1},\phi^{n},\phi^{n-1}}\frac{\phi^{n+1}-\phi^{n}}{\Delta t},\\
\chi^{q^{n},q^{n-1},\phi^{n},\phi^{n-1}}&=&\displaystyle\frac{3q^n-q^{n-1}}{2F(\frac32\phi^n-\frac12\phi^{n-1})+C}F'(\frac32\phi^n-\frac12\phi^{n-1}).
   \end{array}
   \right.
\end{equation}

Similarly, $\phi^{n+1}$ and $\eta^{n+1}$ can be solved step-by-step:
\begin{equation}\label{3s-ieq-e4}
\aligned
&\phi^{n+1}=(I-\frac12\Delta t\mathcal{G}\mathcal{L})^{-1}\phi^n+\frac12(I-\frac12\Delta t\mathcal{G}\mathcal{L})^{-1}\mathcal{G}\mathcal{L}\phi^n+\Delta t(I-\Delta t\mathcal{G}\mathcal{L})^{-1}\mathcal{G}\chi^{q^{n},q^{n-1},\phi^{n},\phi^{n-1}},\\
&\eta^{n+1}=\displaystyle\eta^n+\chi^{q^{n},q^{n-1},\phi^{n},\phi^{n-1}}(\phi^{n+1}-\phi^{n}).
\endaligned
\end{equation}

\begin{remark}\label{3s-ieq-re1}
The above 3S-IEQ scheme can not only be implemented step-by-step, but also has a very obvious advantage than traditional IEQ approach. For solving the discrete scheme based on IEQ approach in \eqref{ieq-e2}, although one only needs to solve a linear system at each time step, the linear system usually involves variable coefficients which change at each time step. However, from \eqref{3s-ieq-e4}, one can see that the discrete scheme based on 3S-IEQ approach leads to linear equations with constant coefficients so it is remarkably easy to implement.
\end{remark}

Taking the inner products of the above equations in \eqref{3s-ieq-e5} with $\mu^{n+\frac12}$, $(\phi^{n+1}-\phi^{n})/\Delta t$ and $1$, we can easily to obtain the following theorem immediately:
\begin{theorem}\label{3s-ieq-th1}
The scheme \eqref{3s-ieq-e5} for the equivalent gradient flow system \eqref{3s-ieq-e3} is second-order accurate, unconditionally energy stable in the sense that
\begin{equation*}
\aligned
\frac{1}{\Delta t}\left[E_{3S-IEQ/CN}^{n+1}-E^{n}_{3S-IEQ/CN}\right]\leq(\mathcal{G}\mu^{n+\frac12},\mu^{n+\frac12})\leq0,
\endaligned
\end{equation*}
where the modified discrete version of the energy is defined by
\begin{equation*}
\aligned
E_{3S-IEQ/CN}^{n}=\frac12(\phi^n,\mathcal{L}\phi^{n})+\int_{\Omega}q^nd\textbf{x}.
\endaligned
\end{equation*}
\end{theorem}
\section{Examples and discussion}
In this section, we give several numerical examples to demonstrate the accuracy, energy stability and efficiency of the proposed 3S-SAV and 3S-IEQ schemes when applying to the some classical gradient flows such as Allen-Cahn equation , Cahn-Hilliard equation, phase field crystal equation and so on. Two comparative studies of traditional SAV/IEQ and 3S-SAV/3S-IEQ approaches are considered to show the accuracy and efficiency. In all examples, we consider the periodic boundary conditions and use a Fourier spectral method in space. To test the efficiency of fast calculation, all the solvers are implemented using Matlab and all the numerical experiments are performed on a computer with 8-GB memory.

\subsection{Allen-Cahn and Cahn-Hilliard equations}
Both Allen-Cahn and Cahn-Hilliard equations are very classical phase field models and have been widely used in many fields involving physics, materials science, finance and image processing \cite{chen2018accurate,chen2018power,du2018stabilized}.

Consider the following Lyapunov energy functional:
\begin{equation}\label{section5_energy1}
E(\phi)=\int_{\Omega}(\frac{\epsilon^2}{2}|\nabla \phi|^2+F(\phi))d\textbf{x},
\end{equation}
where the most commonly used form Ginzburg-Landau double-well type potential is defined as $F(\phi)=\frac{1}{4}(\phi^2-1)^2$.

By applying the variational approach for the free energy \eqref{section5_energy1} leads to
\begin{equation}\label{section5_e_model}
  \left\{
   \begin{array}{rlr}
\displaystyle\frac{\partial \phi}{\partial t}&=M\mathcal{G}\mu,     &(\textbf{x},t)\in\Omega\times J,\\
                                          \mu&=-\epsilon^2\Delta \phi+f(\phi),&(\textbf{x},t)\in\Omega\times J,
   \end{array}
   \right.
  \end{equation}
 where $J=(0,T]$, $M$ is the mobility constant, $\mathcal{G}=-I$ for the Allen-Cahn type system and $\mathcal{G}=\Delta$ for the Cahn-Hilliard type system. $\mu$ is the chemical potential, and $f(\phi)=F^{\prime}(\phi)$.

\textbf{Example 1}: Consider the above Allen-Cahn and Cahn-Hilliard equations in $\Omega=[0,2\pi]$ with $\epsilon=0.1$, $T=0.032$, $M=1$ in Allen-Cahn equation and $M=0.1$ in Cahn-Hilliard equation, and the following initial condition \cite{ShenA}:
\begin{equation*}
\aligned
\phi(x,y,0)=0.05sin(x)sin(y).
\endaligned
\end{equation*}

We use the Fourier spectral Galerkin method for spatial discretization with $N=128$. The true solution is unknown and we therefore use the
Fourier Galerkin approximation in the case $\Delta t=1e-6$ as a reference solution.

For Allen-Cahn equation, we consider first-order time discrete schemes based on both SAV approach in \cite{ShenA} and the proposed 3S-SAV approach in this article. We also give a comparative study of traditional IEQ and 3S-IEQ approaches to show the accuracy and efficiency. For the results of SAV and 3S-SAV schemes, the computational error and convergence rates are shown in Table \ref{tab:tab1}. The numerical results indicate that both SAV and 3S-SAV scheme are indeed of first order in time. However, the CPU time in Table \ref{tab:tab1} shows that the 3S-SAV scheme is about half as time-consuming as SAV scheme. The results of error and convergence rates for IEQ and 3S-IEQ approaches are shown in Table \ref{tab:tab2}. One can see that both accuracy and computational time for 3S-IEQ approach is better than that for traditional IEQ approach. For Cahn-Hilliard model, a comparative study of traditional SAV and 3S-SAV approaches based on second-order Crank-Niclosion scheme is considered and the relative results are shown in Table \ref{tab:tab2}. Similar study for IEQ and 3S-IEQ approaches is shown in Table \ref{tab:tab4}. One can see that these four methods obtain almost identical error and convergence rates. Similar as before, the 3S-SAV scheme also saves half the time compared with SAV scheme, which is essentially same for IEQ and 3S-IEQ approaches.
\begin{table}[h!b!p!]
\small
\centering
\caption{\small The $L_2$ errors, convergence rates for first-order scheme in time for SAV and 3S-SAV approaches of Allen-Cahn equation.}\label{tab:tab1}
\begin{tabular}{cccccccccc}
\hline
          &&SAV&&&&3S-SAV&\\
\cline{1-8}
$\Delta t$          &$L_2$ error&Rate&Cpu-Time(s)&&$L_2$ error&Rate&Cpu-Time(s)\\
\cline{1-8}
$1.6e-4$ &2.4812e-5   &---   &1.14    &&2.4813e-5   &---   &0.67\\
$8e-5$   &1.2328e-5   &1.0091&2.09    &&1.2328e-5   &1.0091&1.21\\
$4e-5$   &6.0860e-6   &1.0183&4.31    &&6.0862e-6   &1.0183&2.23\\
$2e-5$   &2.9649e-6   &1.0375&8.40    &&2.9651e-6   &1.0375&4.53\\
$1e-5$   &1.4044e-6   &1.0780&17.89   &&1.4045e-6   &1.0780&9.01\\
\cline{1-8}
\end{tabular}
\end{table}

\begin{table}[h!b!p!]
\small
\centering
\caption{\small The $L_2$ errors, convergence rates for second-order scheme in time for SAV and 3S-SAV approaches of Cahn-Hilliard equation.}\label{tab:tab2}
\begin{tabular}{cccccccccc}
\hline
          &&SAV&&&&3S-SAV&\\
\cline{1-8}
$\Delta t$          &$L_2$ error&Rate&Cpu-Time(s)&&$L_2$ error&Rate&Cpu-Time(s)\\
\cline{1-8}
$1.6e-4$ &3.9566e-9   &---    &1.36    &&3.9566e-9   &---   &0.73\\
$8e-5$   &9.8907e-10   &2.0001&2.61    &&9.8907e-10  &2.0001&1.48\\
$4e-5$   &2.4714e-10   &2.0007&5.10    &&2.4714e-10  &2.0007&2.99\\
$2e-5$   &6.1646e-11  &2.0032 &10.06   &&6.1649e-11  &2.0032&5.72\\
$1e-5$   &1.5394e-11  &2.0116 &20.20   &&1.5400e-11  &2.0016&11.76\\
\cline{1-8}
\end{tabular}
\end{table}

\begin{table}[h!b!p!]
\small
\centering
\caption{\small The $L_2$ errors, convergence rates for first-order scheme in time for IEQ and 3S-SAV approaches of Allen-Cahn equation.}\label{tab:tab3}
\begin{tabular}{cccccccccc}
\hline
          &&IEQ&&&&3S-IEQ&\\
\cline{1-8}
$\Delta t$          &$L_2$ error&Rate&Cpu-Time(s)&&$L_2$ error&Rate&Cpu-Time(s)\\
\cline{1-8}
$1.6e-4$ &2.5010e-5   &---   &0.68    &&2.4813e-5   &---   &0.37\\
$8e-5$   &1.2534e-5   &0.9966&1.32    &&1.2328e-5   &1.0091&0.73\\
$4e-5$   &6.3095e-6   &0.9902&2.55    &&6.0862e-6   &1.0183&1.56\\
$2e-5$   &3.2235e-6   &0.9688&5.21    &&2.9651e-6   &1.0375&3.19\\
$1e-5$   &1.7284e-6   &0.8992&10.07   &&1.4045e-6   &1.0780&6.09\\
\cline{1-8}
\end{tabular}
\end{table}

\begin{table}[h!b!p!]
\small
\centering
\caption{\small The $L_2$ errors, convergence rates for first-order scheme in time for IEQ and 3S-IEQ approaches of Cahn-Hilliard equation.}\label{tab:tab4}
\begin{tabular}{cccccccccc}
\hline
          &&IEQ&&&&3S-IEQ&\\
\cline{1-8}
$\Delta t$          &$L_2$ error&Rate&Cpu-Time(s)&&$L_2$ error&Rate&Cpu-Time(s)\\
\cline{1-8}
$1.6e-4$ &3.9591e-9    &---   &0.59    &&3.9595e-9    &---   &0.35\\
$8e-5$   &9.9003e-10   &1.9996&1.23    &&9.8587e-10   &2.0058&0.78\\
$4e-5$   &2.4763e-10   &1.9992&2.51    &&2.4523e-10   &2.0072&1.59\\
$2e-5$   &6.1976e-11   &1.9984&4.78    &&6.0983e-11   &2.0076&3.08\\
$1e-5$   &1.5658e-11   &1.9848&9.75    &&1.5164e-11   &2.0077&6.08\\
\cline{1-8}
\end{tabular}
\end{table}

\textbf{Example 2}: In the following, we solve a benchmark problem for the Allen-Cahn equation which can be seen in many articles such as \cite{ShenA}. We take $\epsilon=0.01$, $M=1$. The initial condition is chosen as
\begin{equation*}
\aligned
\phi_0(x,y,0)=\sum\limits_{i=1}^2-\tanh\left(\frac{\sqrt{(x-x_i)^2+(y-y_i)^2}-R_i}{\sqrt{2}\epsilon}\right)+1.
\endaligned
\end{equation*}
with the radius $R_1=0.15$, $(x_1,y_1)=(0.35,0.35)$ and $R_2=0.2$, $(x_2,y_2)=(0.6,0.6)$. Initially, two bubbles, centered at $(0.35,0.35)$ and $(0.6,0.6)$, respectively, are osculating.

\begin{figure}[htp]
\centering
\subfigure[t=0]{
\includegraphics[width=3.8cm,height=3.8cm]{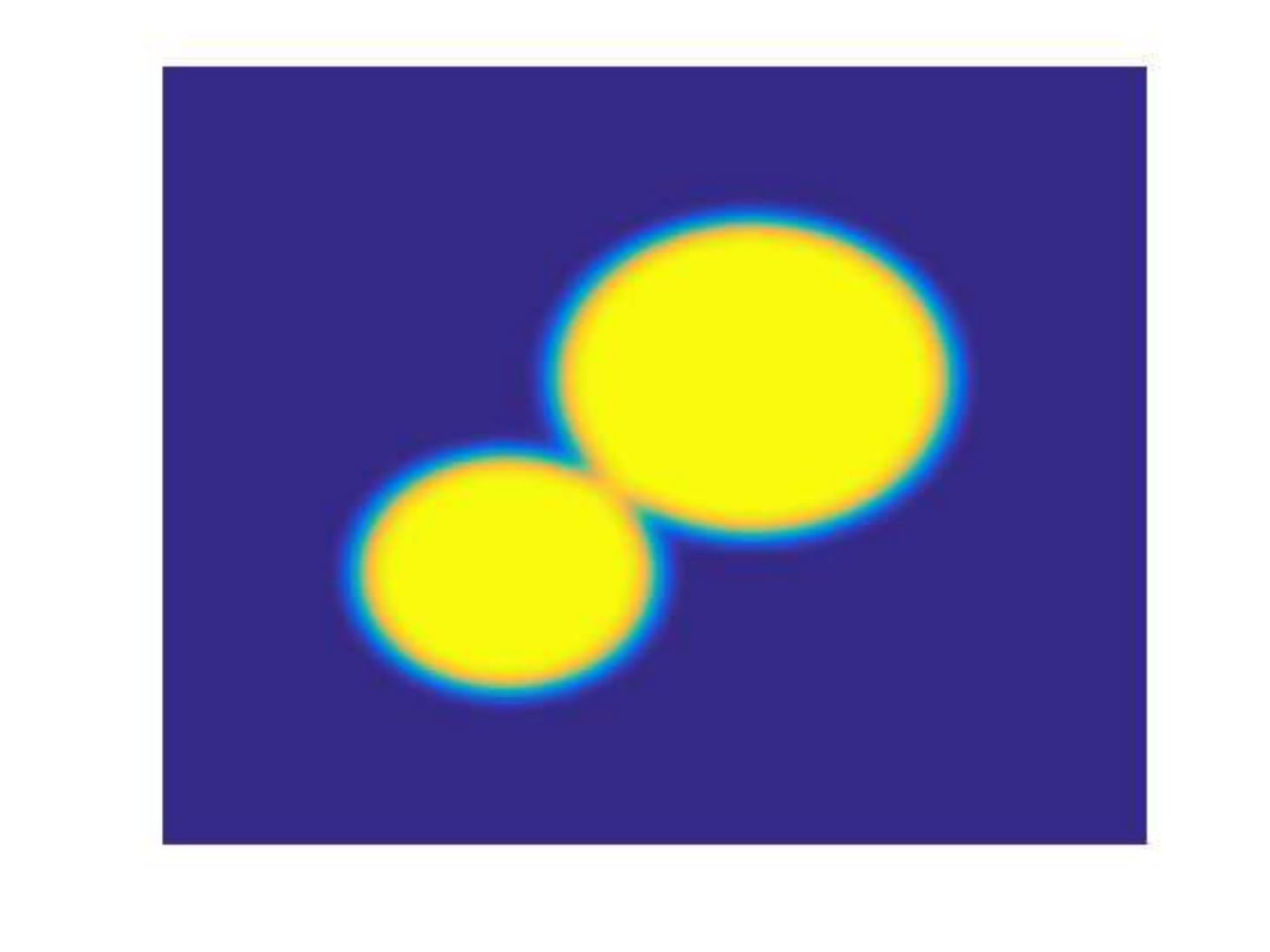}
}
\subfigure[t=2]
{
\includegraphics[width=3.8cm,height=3.8cm]{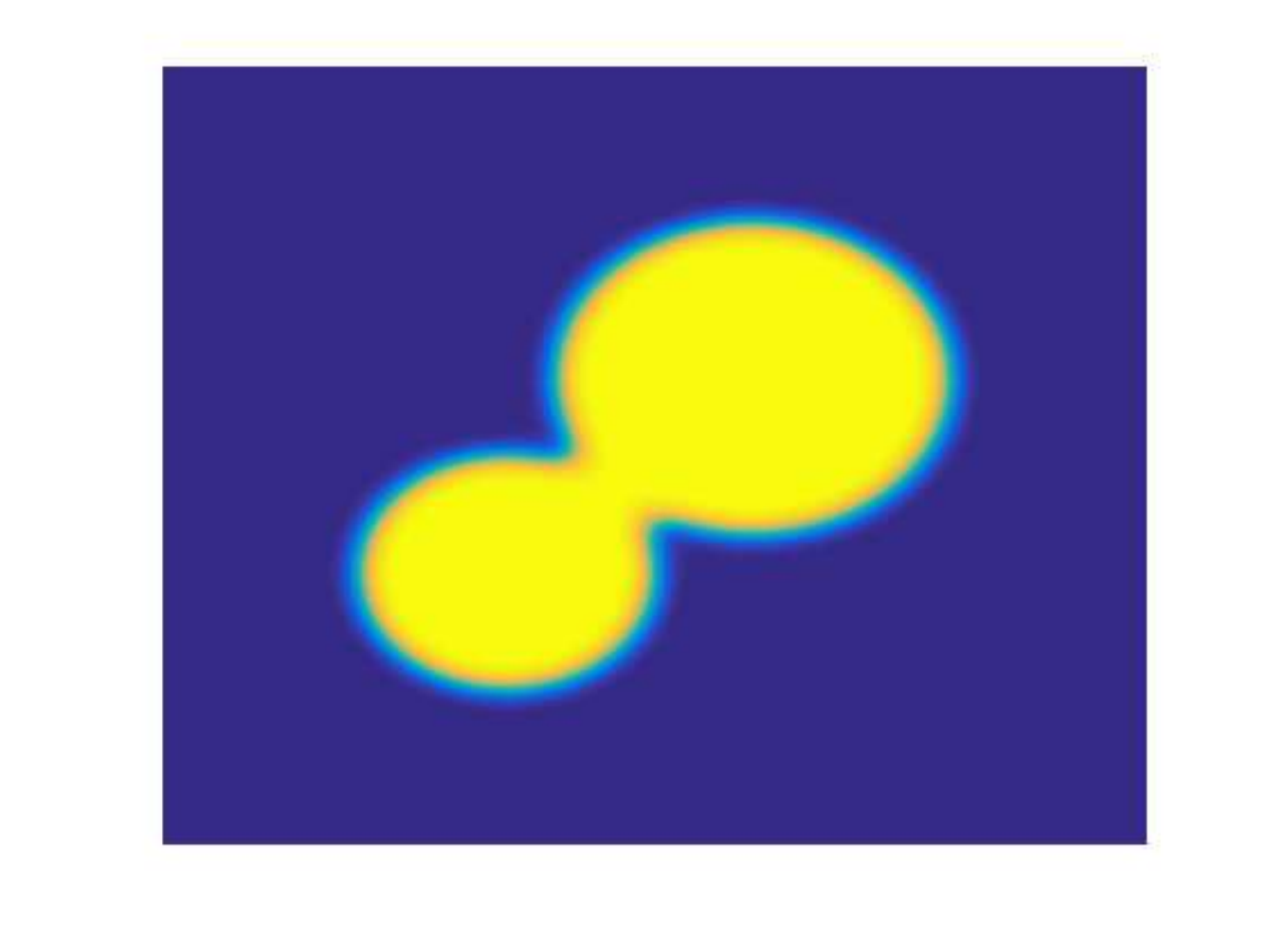}
}
\subfigure[t=5]
{
\includegraphics[width=3.8cm,height=3.8cm]{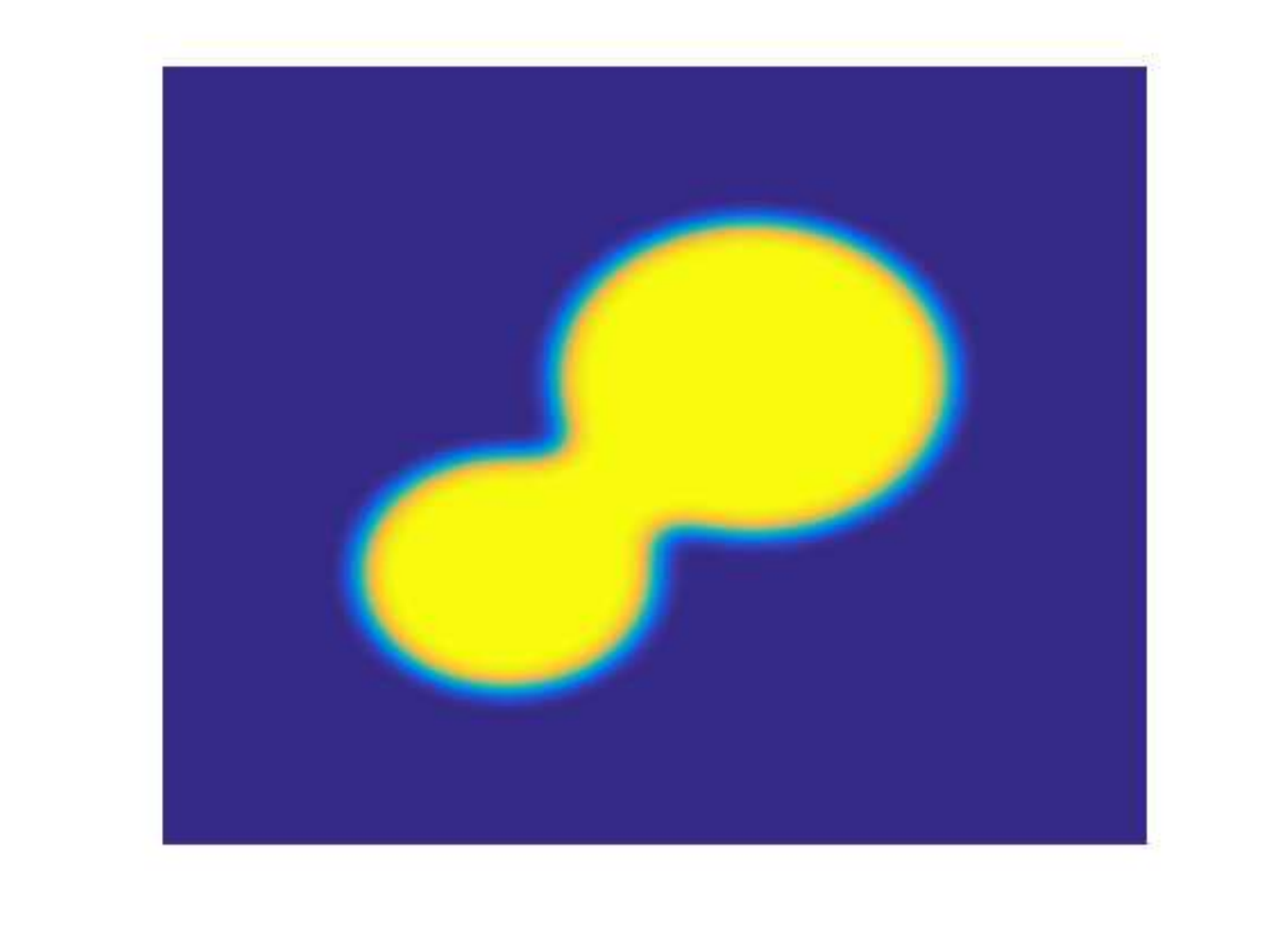}
}
\subfigure[t=50]
{
\includegraphics[width=3.8cm,height=3.8cm]{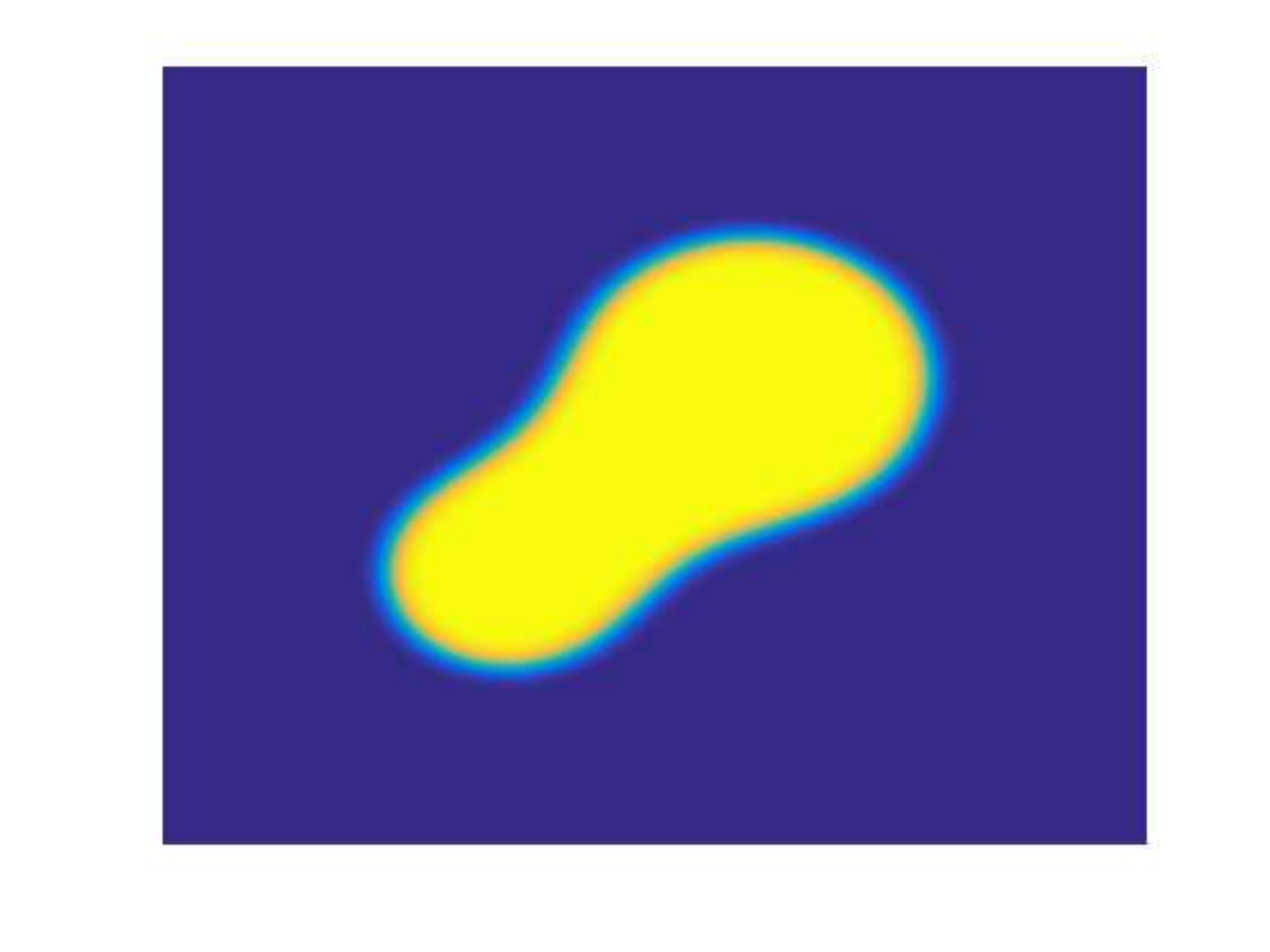}
}
\quad
\subfigure[t=100]
{
\includegraphics[width=3.8cm,height=3.8cm]{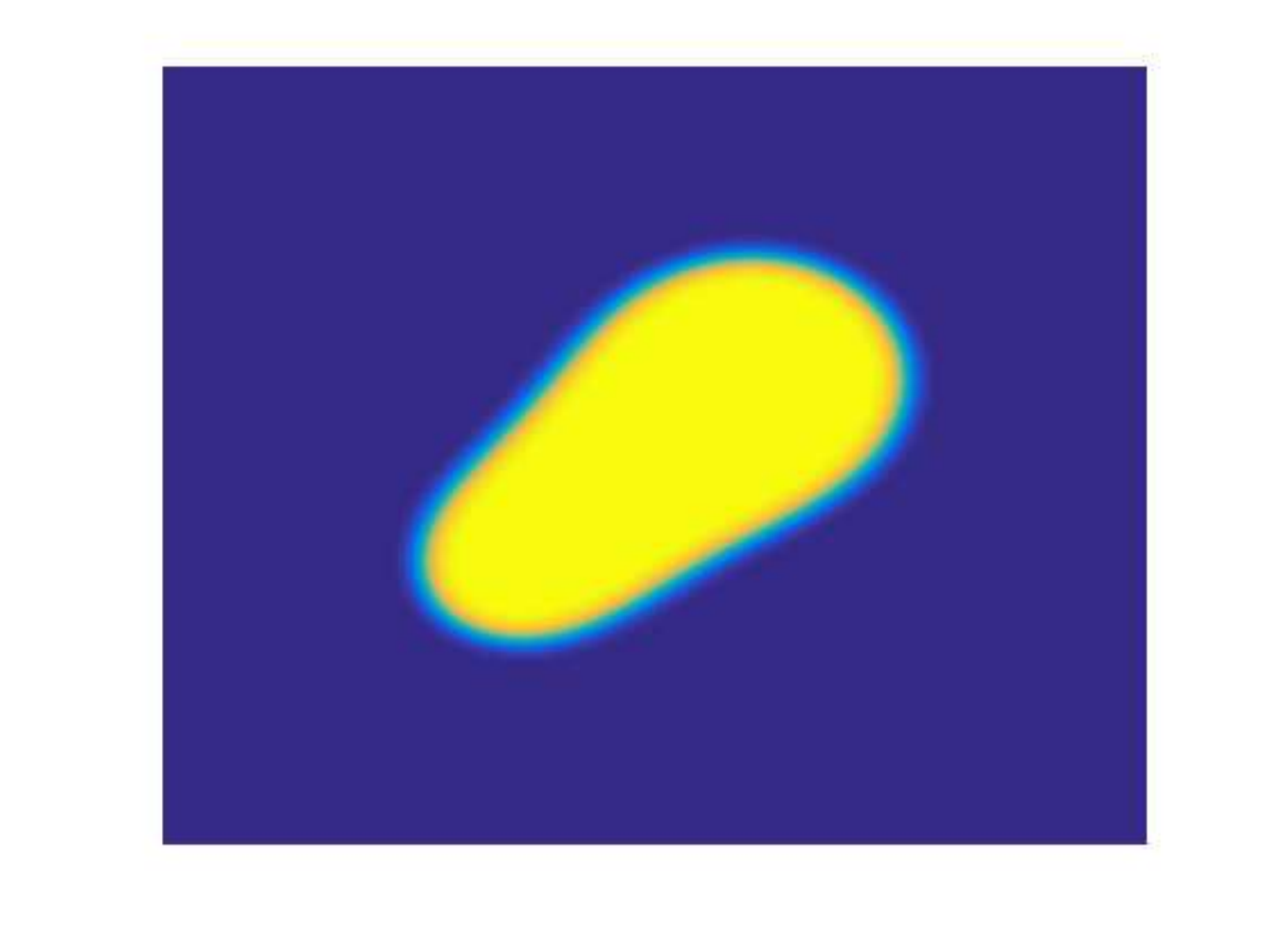}
}
\subfigure[t=250]{
\includegraphics[width=3.8cm,height=3.8cm]{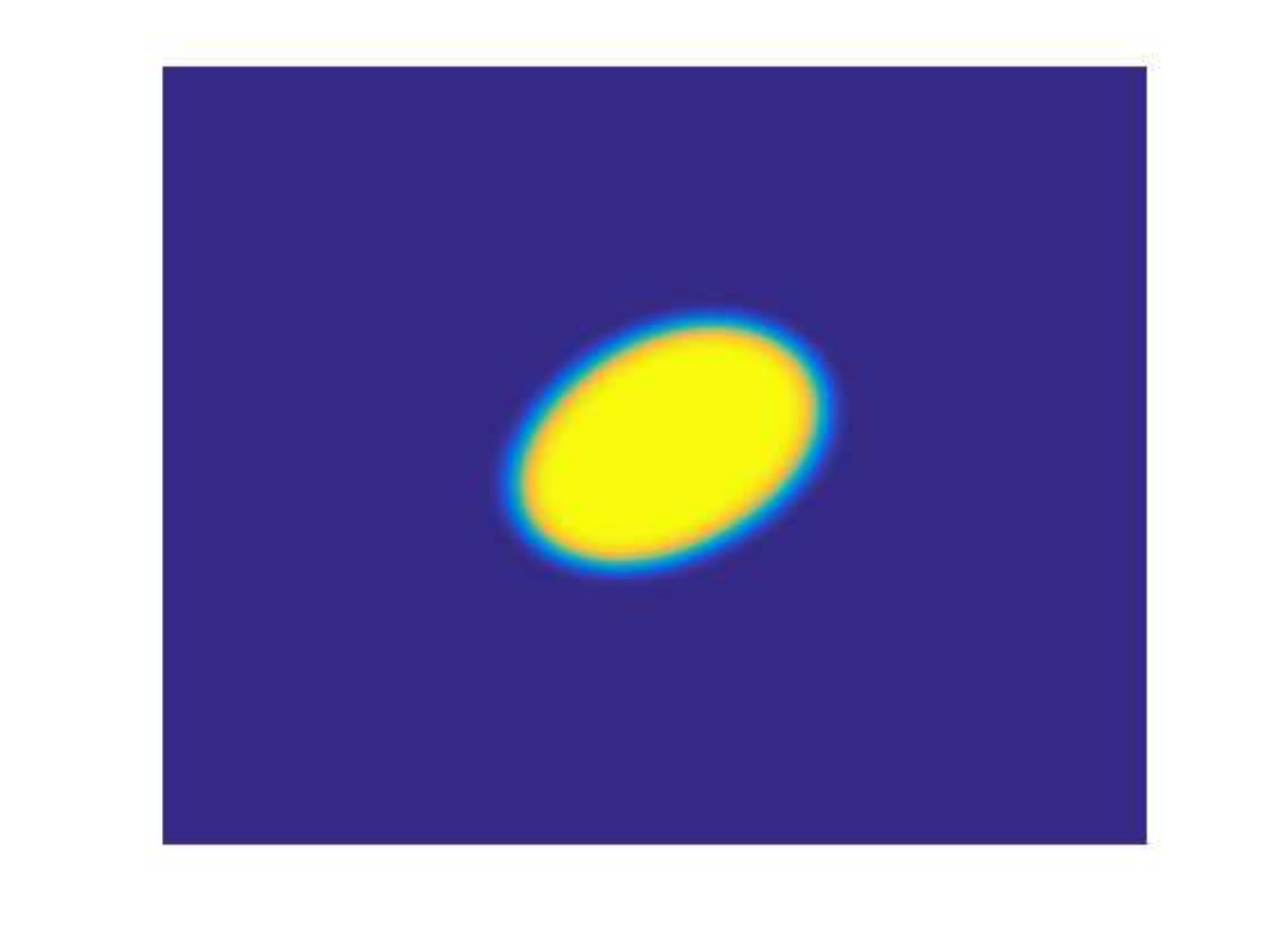}
}
\subfigure[t=340]
{
\includegraphics[width=3.8cm,height=3.8cm]{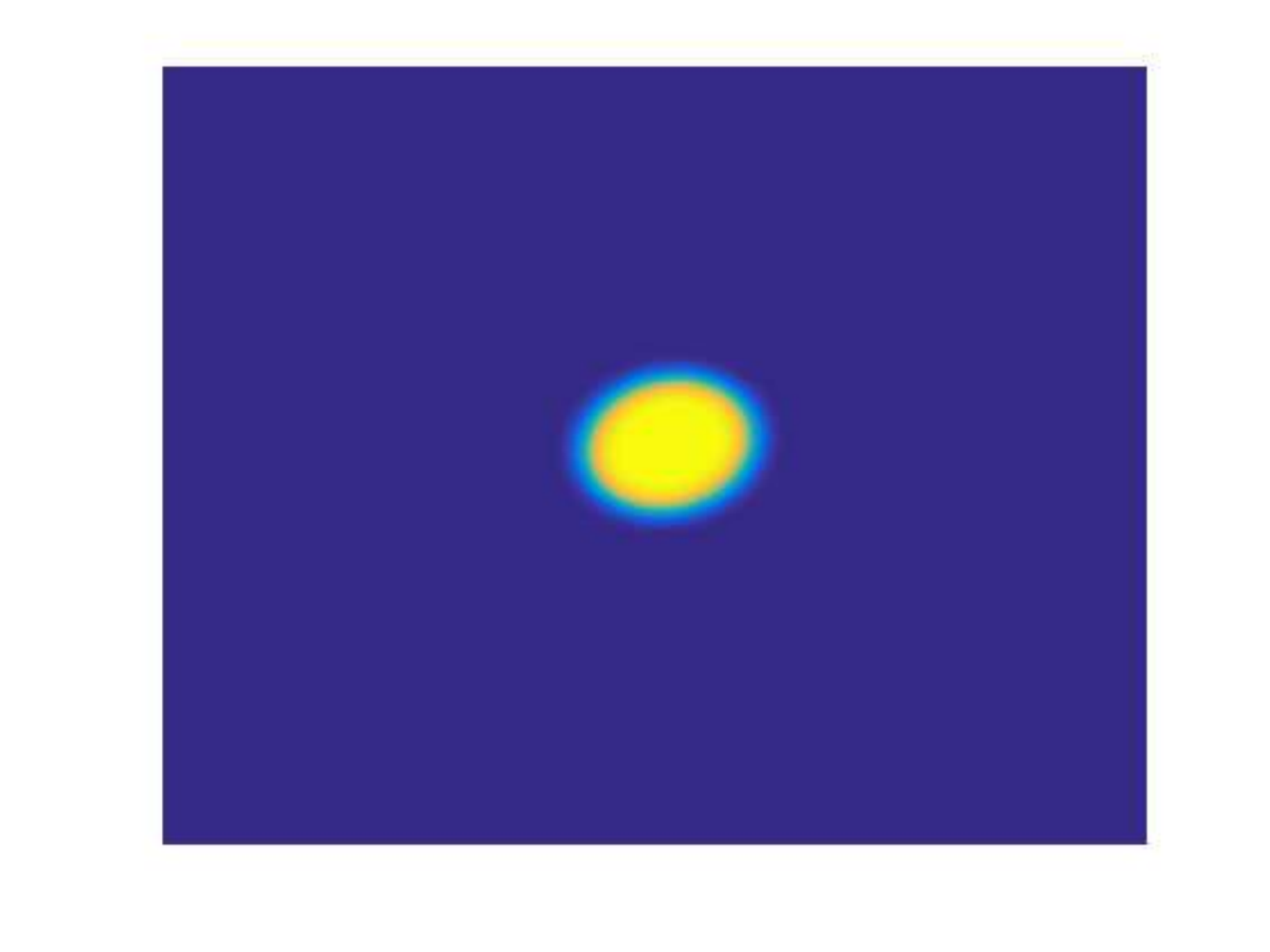}
}
\subfigure[t=380]
{
\includegraphics[width=3.8cm,height=3.8cm]{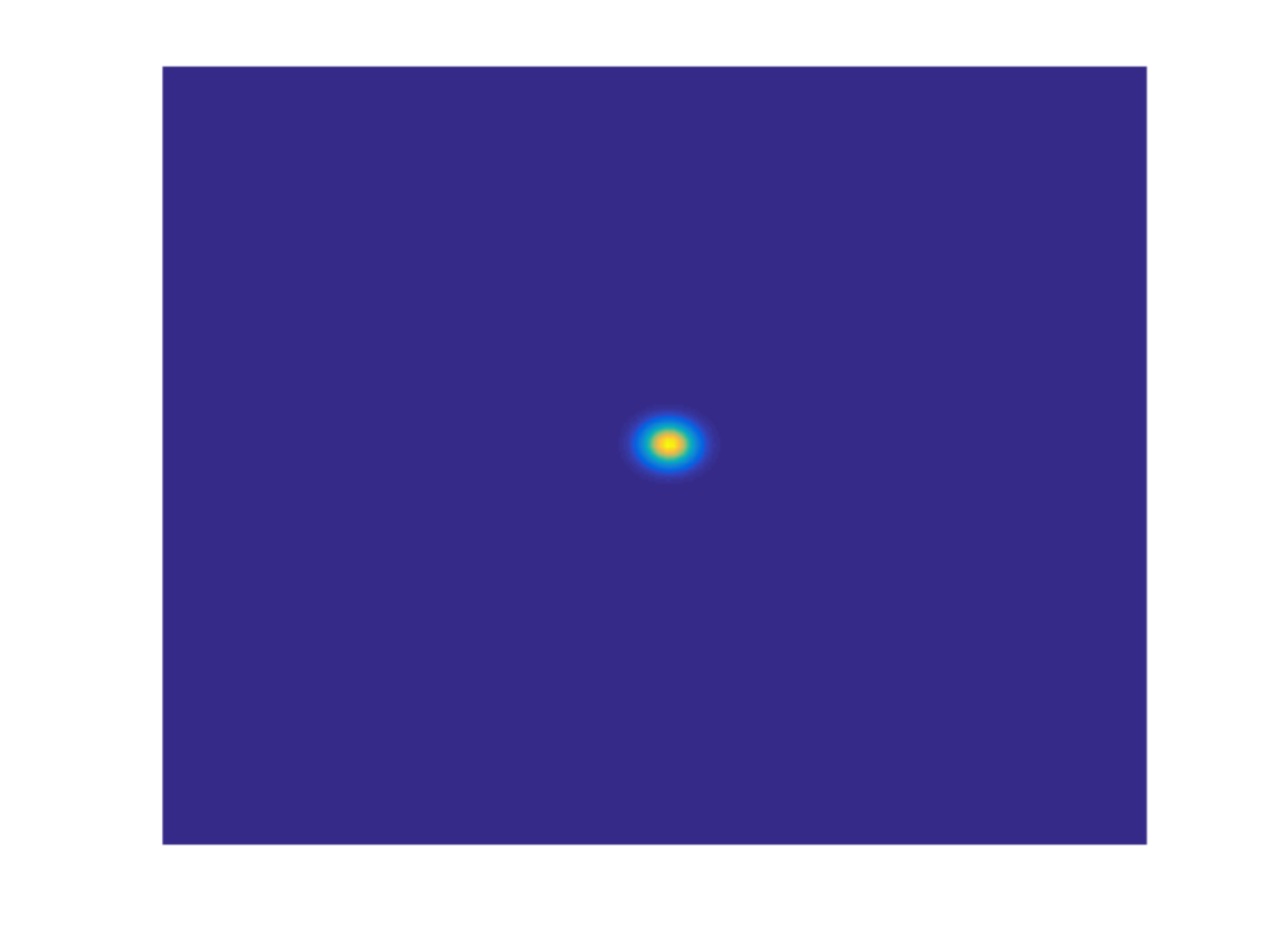}
}
\caption{Snapshots of the phase variable $\phi$ are taken at t=0, 2, 5, 50, 100, 250, 340 and 380 with $\Delta t=0.1$ for example 2.}\label{fig:fig1}
\end{figure}

\begin{figure}[htp]
\centering
\includegraphics[width=7cm,height=7cm]{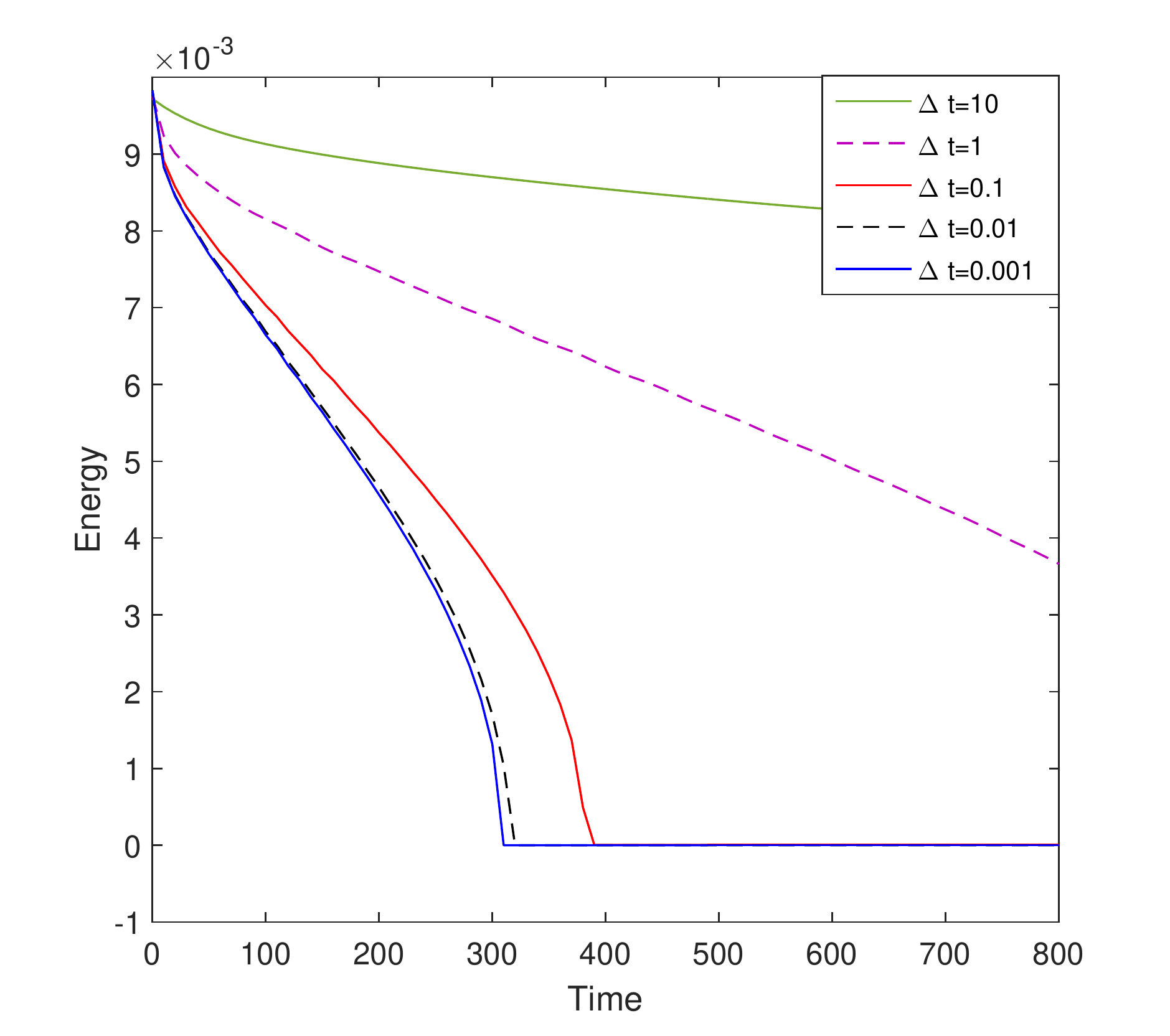}
\includegraphics[width=7cm,height=7cm]{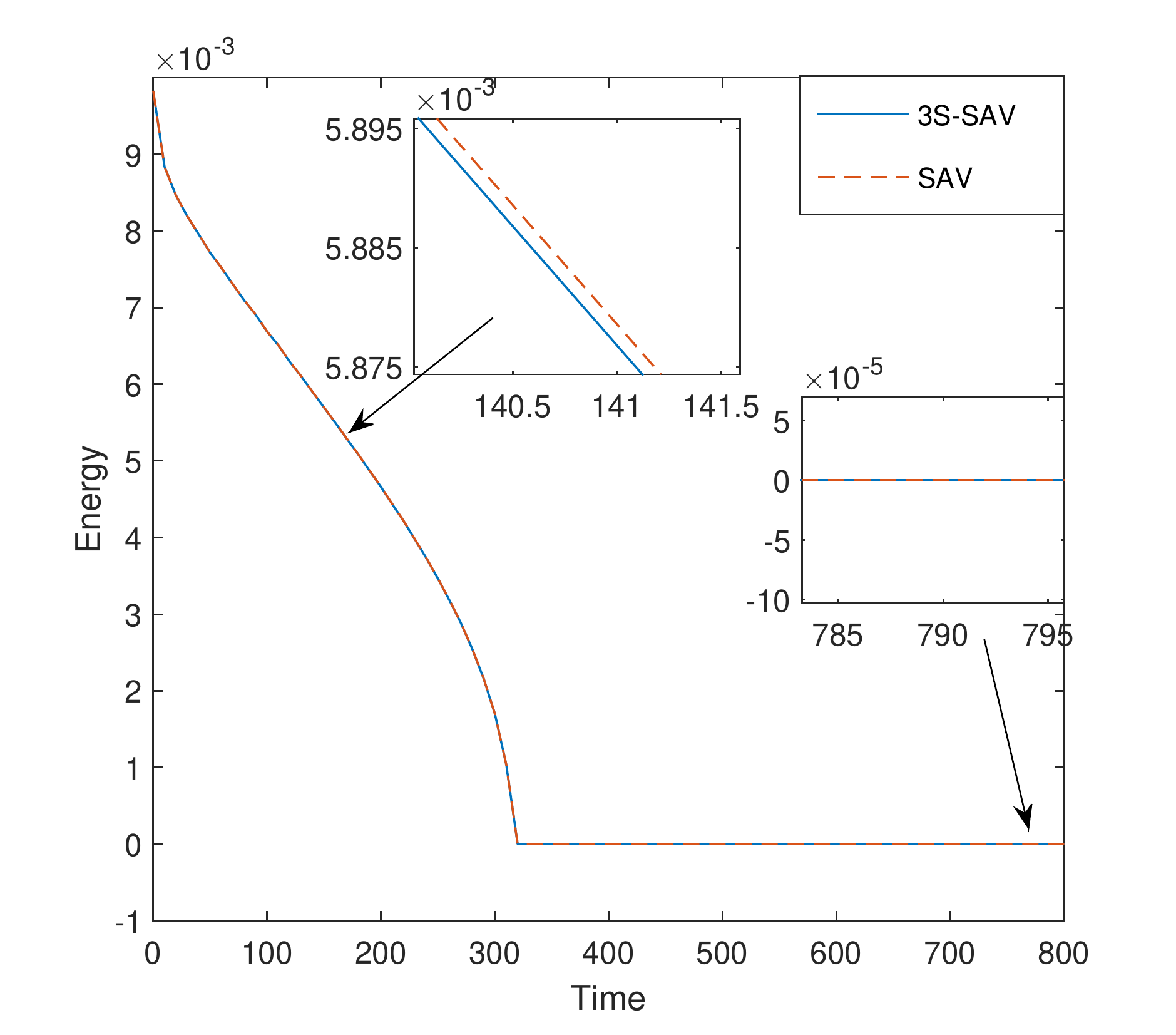}
\caption{Left: time evolution of the energy functional for five different time steps of $\Delta t=0.001$, $0.01$, $0.1$, $1$ and $10$. Right: energy evolution of 3S-SAV and SAV approaches for example 2 with $\Delta t=0.01$.}\label{fig:fig2}
\end{figure}

In Figure \ref{fig:fig1}, It can be found that as time evolves, the two bubbles coalesce into a single bubble, then, shrinks and finally disappears. This process indicates than the Allen-Cahn equation does not conserve mass. In Figure \ref{fig:fig2}, we plot the time evolution of the energy functional with different time step size of $\Delta t=0.001$, $0.01$, $0.1$, $1$ and $10$ by using the first order scheme based on 3S-SAV approach. All energy curves show the monotonic decays for all time steps that confirms that the algorithm 3S-SAV is unconditionally energy stable. Time evolution of the total free energy based on SAV and 3s-SAV approaches is computed by using the time step $\Delta t=0.01$ in the right figure in Figure \ref{fig:fig2}. These two energy curves are almost identical and they both decay monotonically at all times.

\textbf{Example 3}: In the following, we solve a benchmark problem for the Cahn-Hilliard equation on $[0,2\pi)^2$ which can also be seen in many articles such as \cite{shen2018scalar}. When using SAV approach to simulate Cahn-Hilliard model, we need to specify the operators $\mathcal{L}=-\epsilon^2\Delta+\beta$ and $F(\phi)=\frac{1}{4}(\phi^2-1-\beta)^2$ to obtain stable simulation \cite{ShenA}. We take $\epsilon=0.04$, $M=0.1$, $\beta=4$ and discretize the space by the Fourier spectral method with $128\times128$ modes. The initial condition is chosen as the following
\begin{equation*}
\aligned
\phi_0(x,y,0)=0.25+0.4Rand(x,y),
\endaligned
\end{equation*}
where $Rand(x,y)$ is a randomly generated function.

\begin{figure}[htp]
\centering
\subfigure[t=2]{
\includegraphics[width=3.8cm,height=3.8cm]{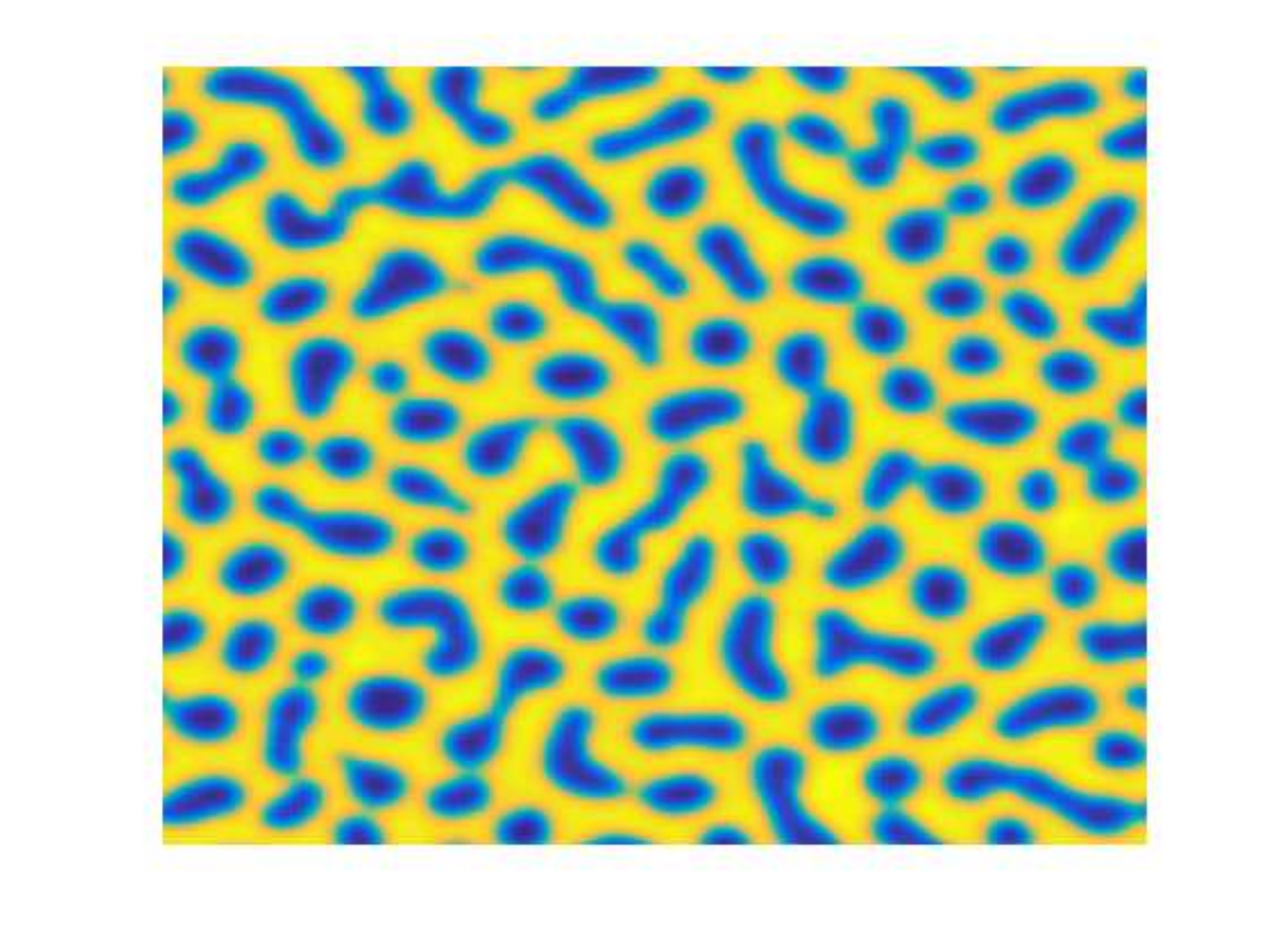}
}
\subfigure[t=10]
{
\includegraphics[width=3.8cm,height=3.8cm]{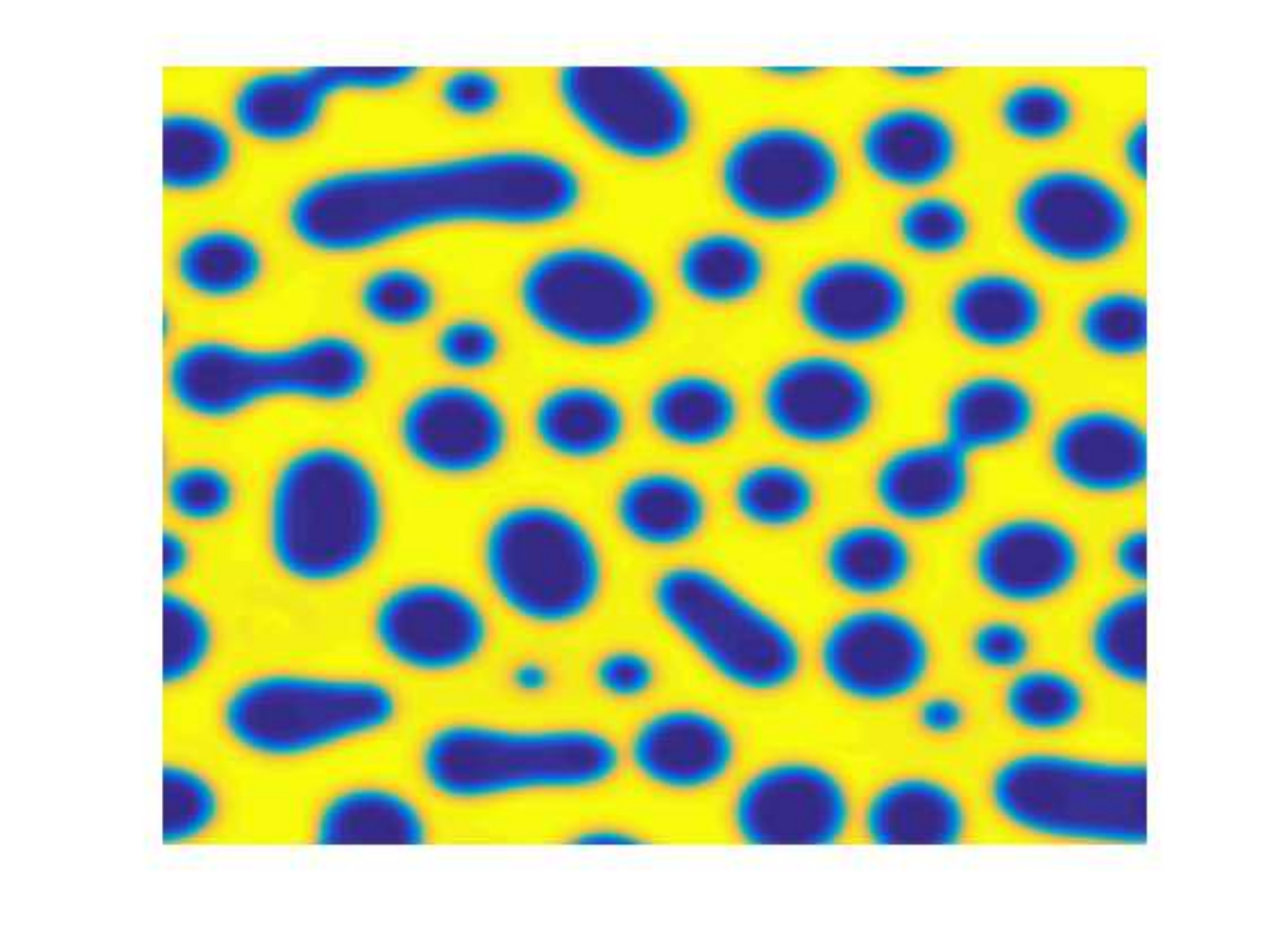}
}
\subfigure[t=40]
{
\includegraphics[width=3.8cm,height=3.8cm]{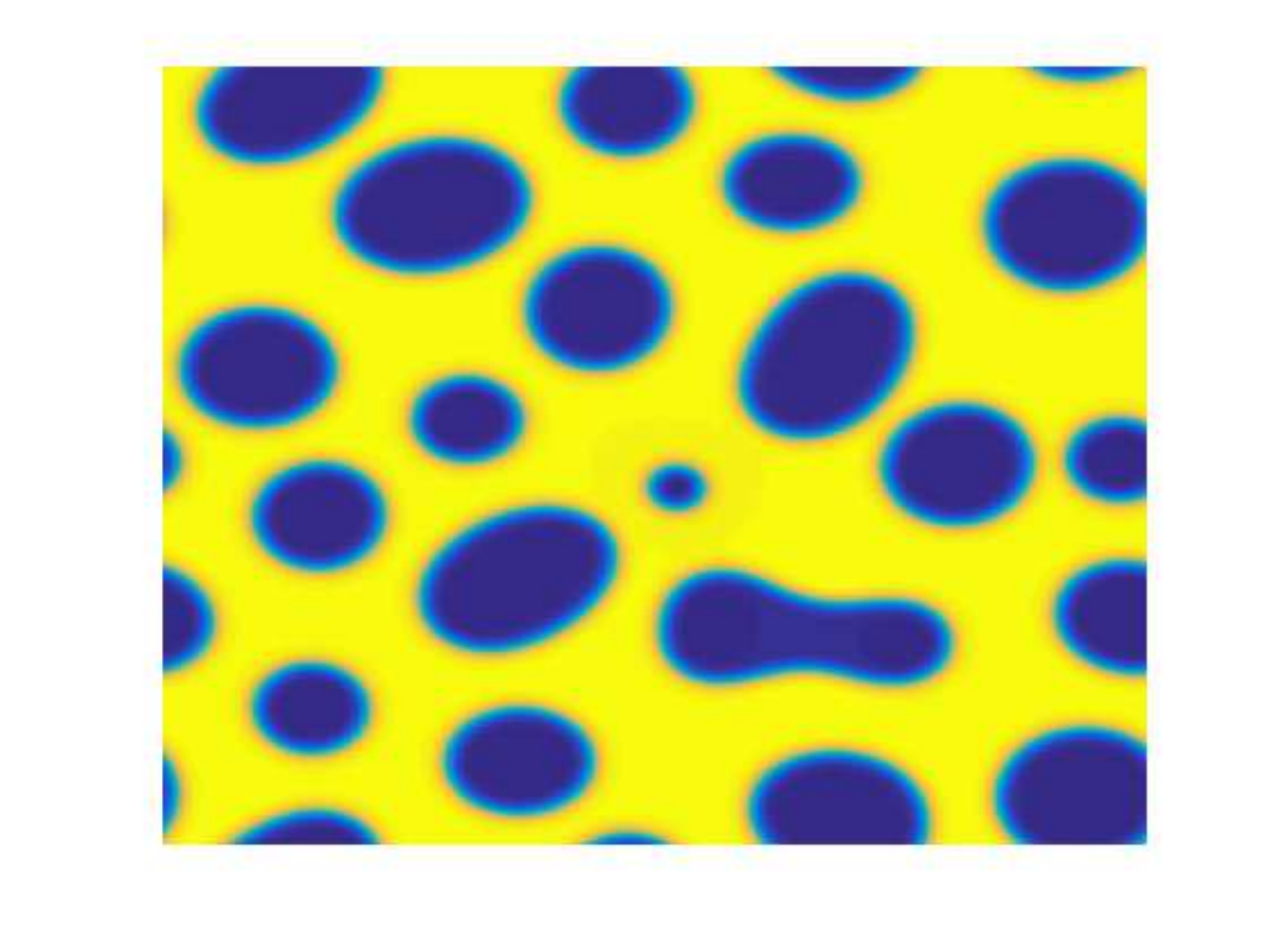}
}
\subfigure[t=100]
{
\includegraphics[width=3.8cm,height=3.8cm]{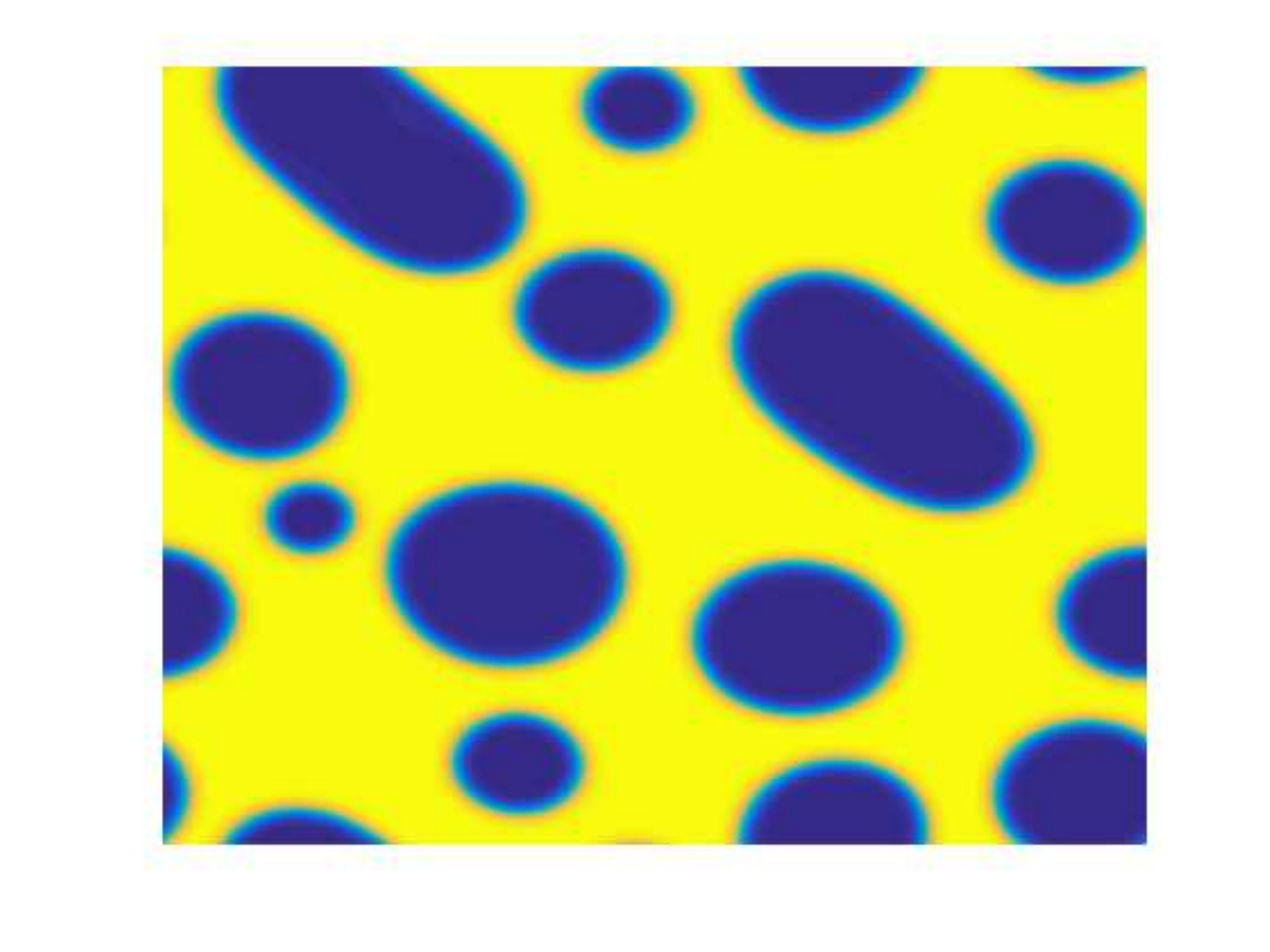}
}
\quad
\subfigure[t=200]
{
\includegraphics[width=3.8cm,height=3.8cm]{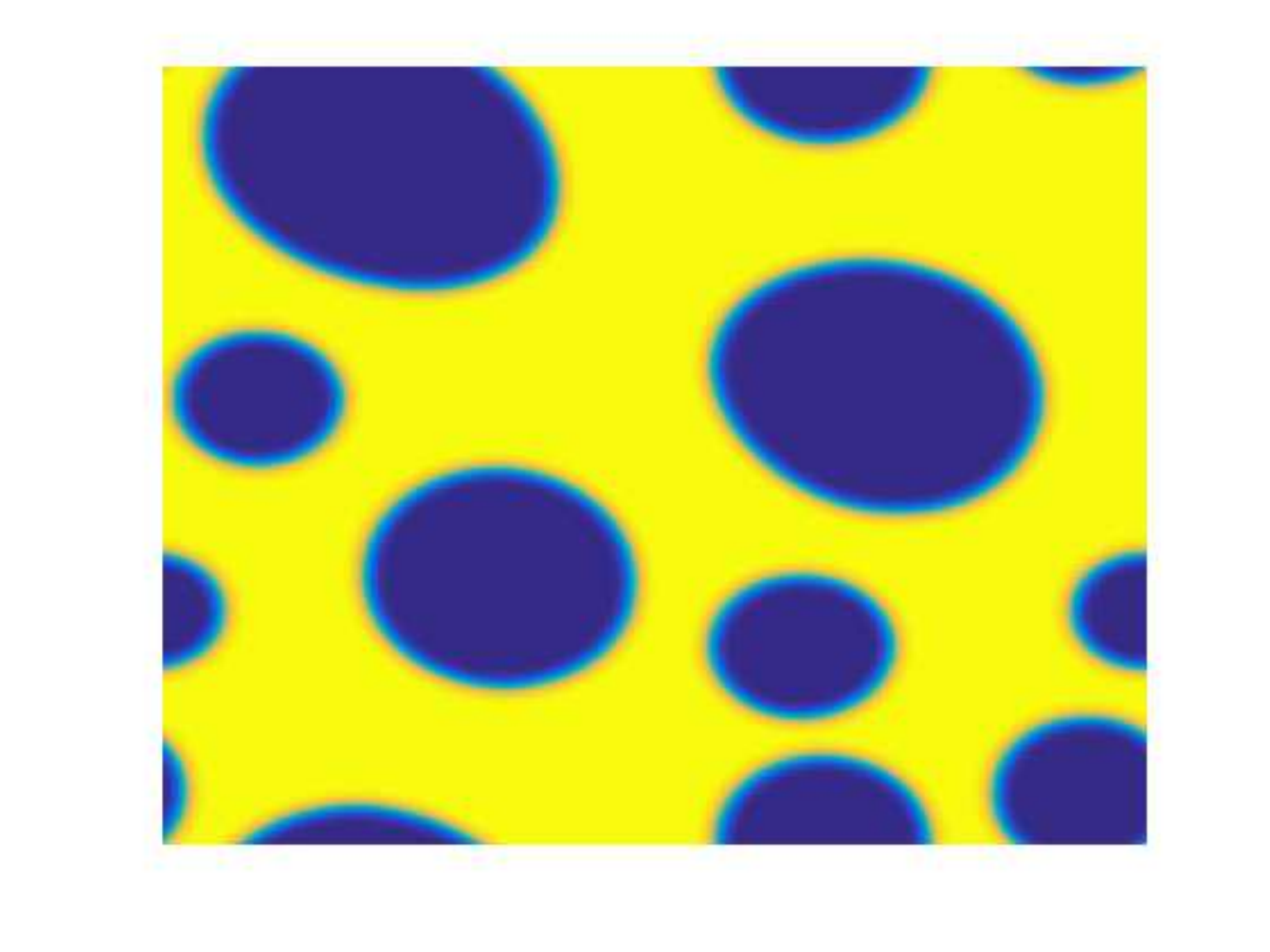}
}
\subfigure[t=400]{
\includegraphics[width=3.8cm,height=3.8cm]{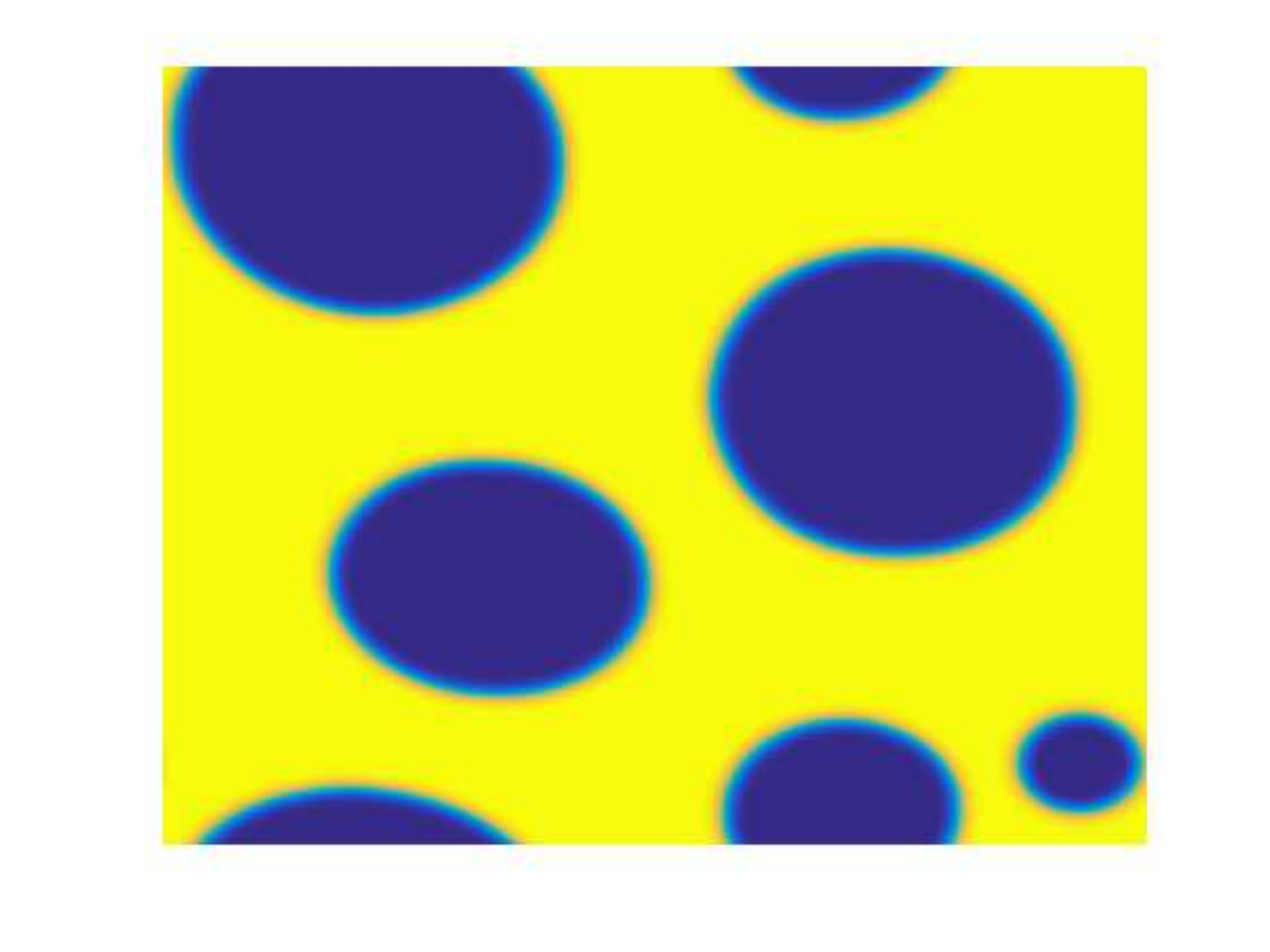}
}
\subfigure[t=1500]
{
\includegraphics[width=3.8cm,height=3.8cm]{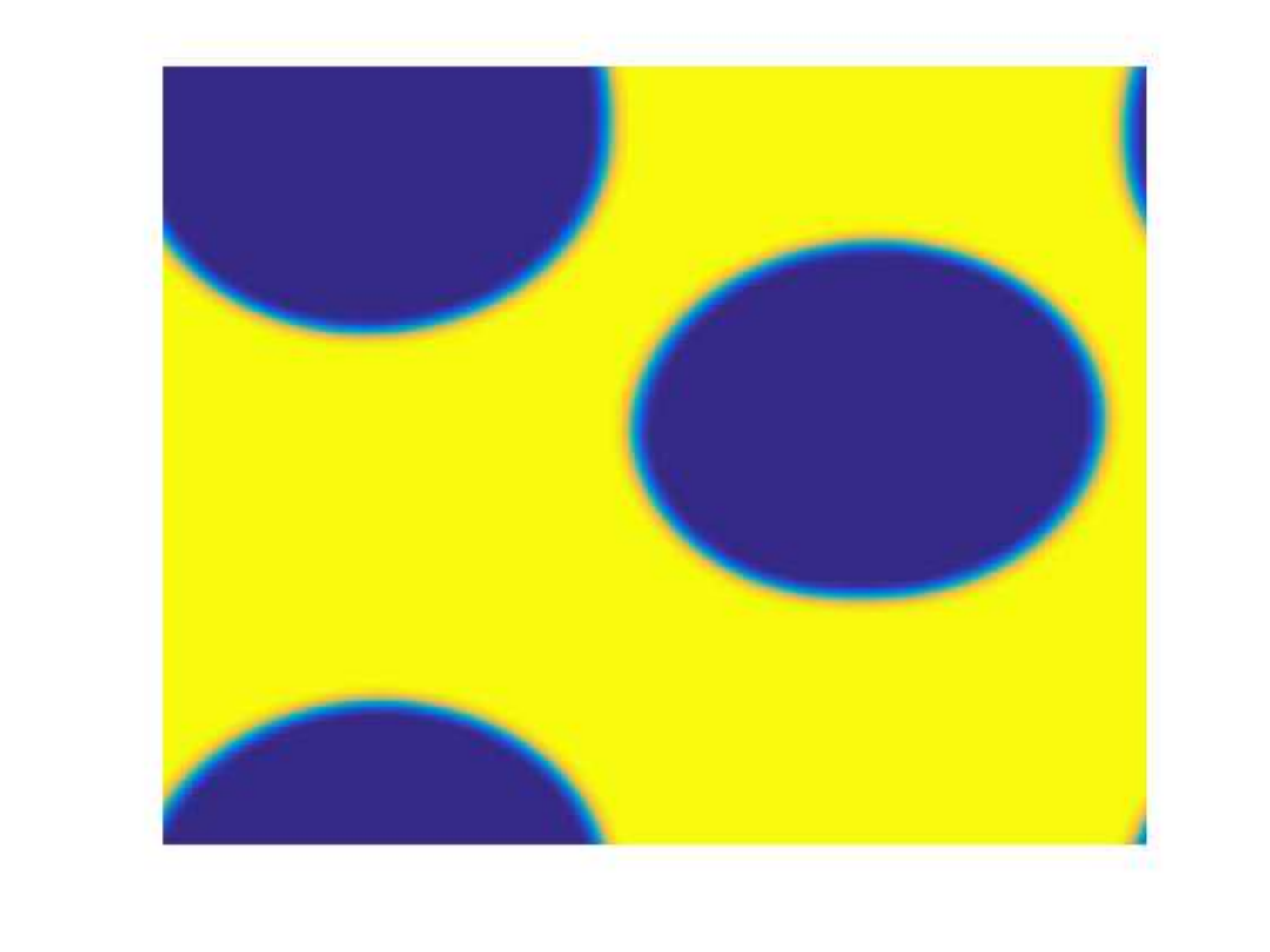}
}
\subfigure[t=8000]
{
\includegraphics[width=3.8cm,height=3.8cm]{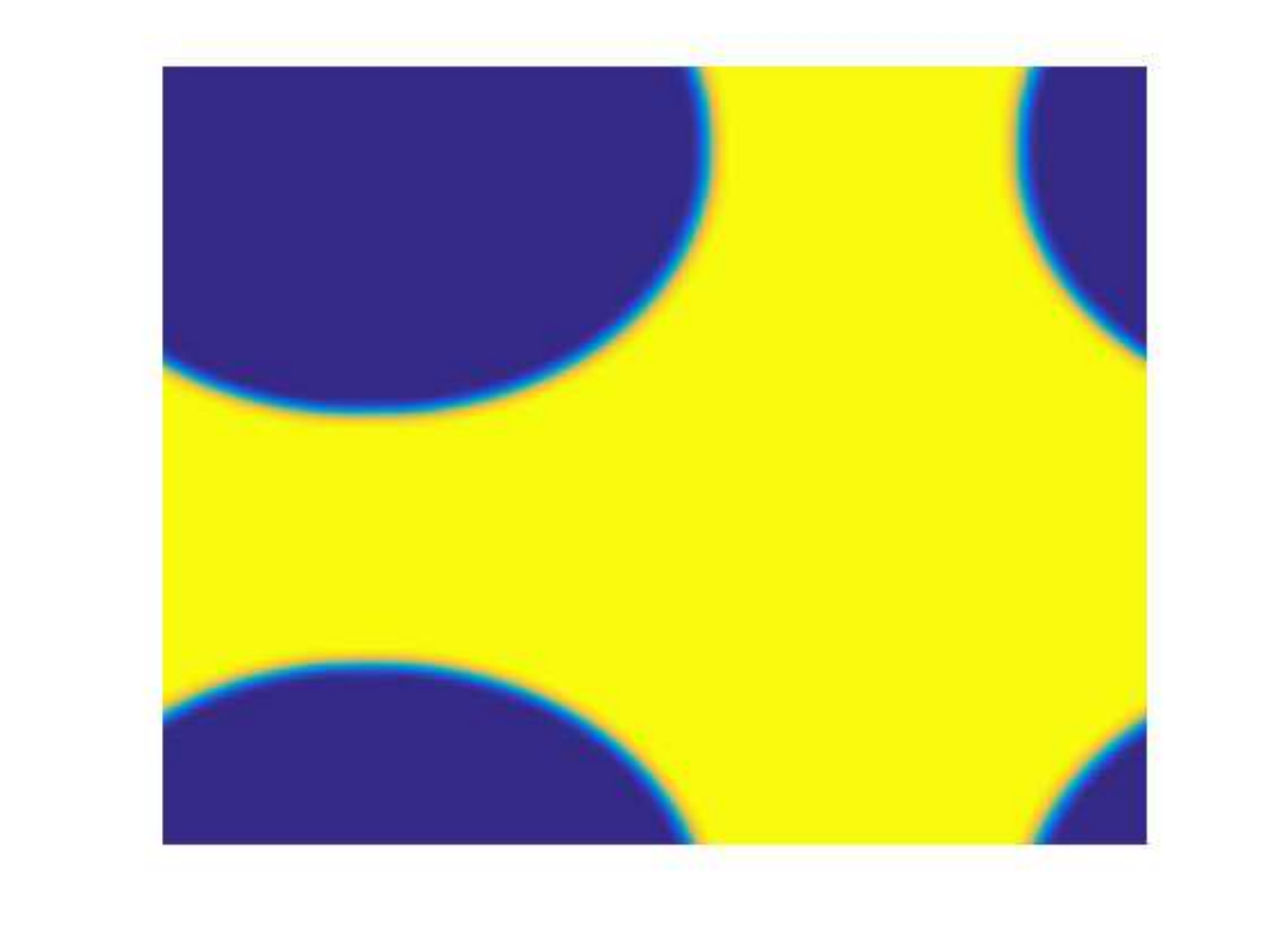}
}
\caption{Snapshots of the phase variable $\phi$ are taken at t=2, 10, 40, 100, 200, 400, 1500 and 8000 with $\Delta t=0.1$ for example 3.}\label{fig:fig3}
\end{figure}

Snapshots of the phase variable $\phi$ taken at $t=2$, $10$, $40$, $100$, $200$, $400$, $1500$ and $8000$ with $\Delta t=0.1$ are shown in Figure \ref{fig:fig3}. The phase separation and coarsening process can be observed very simply which is consistent with the results in \cite{ShenA}.
\subsection{Phase field crystal equations}
Elder \cite{elder2002modeling} firstly proposed the phase field crystal (PFC) model based on density functional theory in 2002. This model can simulate the evolution of crystalline microstructure on atomistic length and diffusive time scales. It naturally incorporates elastic and plastic deformations and multiple crystal orientations, and can be applied to many different physical phenomena.

In particular, consider the following Swift-Hohenberg free energy:
\begin{equation*}
E(\phi)=\int_{\Omega}\left(\frac{1}{4}\phi^4+\frac{1}{2}\phi\left(-\epsilon+(1+\Delta)^2\right)\phi\right)d\textbf{x},
\end{equation*}
where $\textbf{x} \in \Omega \subseteq \mathbb{R}^d$, $\phi$ is the density field and $\epsilon$ is a positive bifurcation constant with physical significance. $\Delta$ is the Laplacian operator.

Considering a gradient flow in $H^{-1}$, one can obtain the phase field crystal equation under the constraint of mass conservation as follows:
\begin{equation*}
\frac{\partial \phi}{\partial t}=\Delta\mu=\Delta\left(\phi^3-\epsilon\phi+(1+\Delta)^2\phi\right), \quad(\textbf{x},t)\in\Omega\times Q,
\end{equation*}
which is a sixth-order nonlinear parabolic equation and can be applied to simulate various phenomena such as crystal growth, material hardness and phase transition. Here $Q=(0,T]$, $\mu=\frac{\delta E}{\delta \phi}$ is called the chemical potential.

Next, we plan to simulate the phase transition behavior of the phase field crystal model. The similar numerical example can be found in many articles such as \cite{li2017efficient,yang2017linearly}.

\textbf{Example 4}: The initial condition is
\begin{equation*}
\aligned
&\phi_0(x,y)=0.07+0.07Rand(x,y),
\endaligned
\end{equation*}
where the $Rand(x,y)$ is the random number in $[-1,1]$ with zero mean. The order parameter is $\epsilon=0.025$, the computational domain $\Omega=[0,128]^2$. we set $256^2$ Fourier modes to discretize the two dimensional space.

\begin{figure}[htp]
\centering
\begin{tabular}{ccccc}
SAV&\includegraphics[width=3.6cm,height=3.6cm]{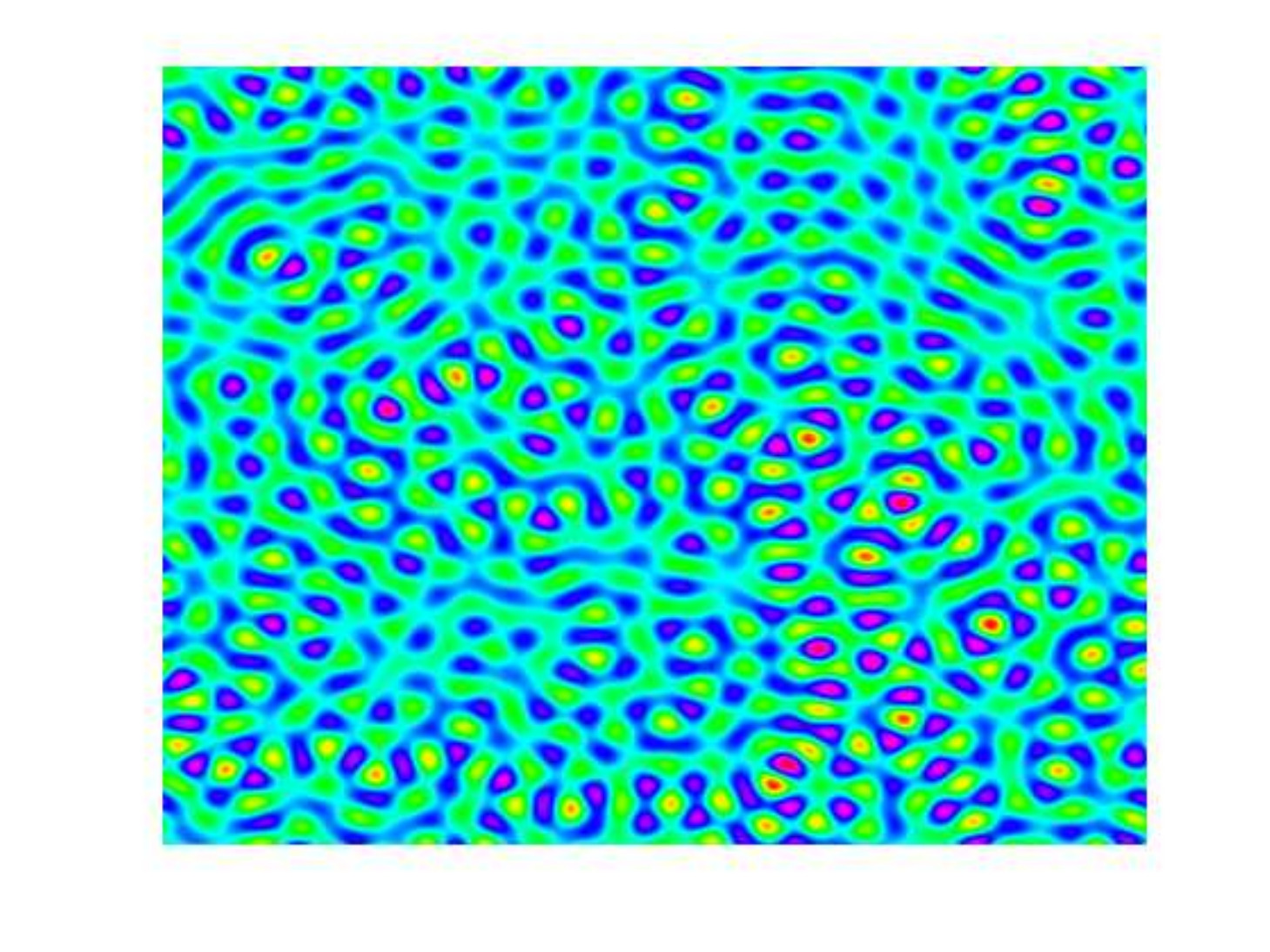}
\includegraphics[width=3.6cm,height=3.6cm]{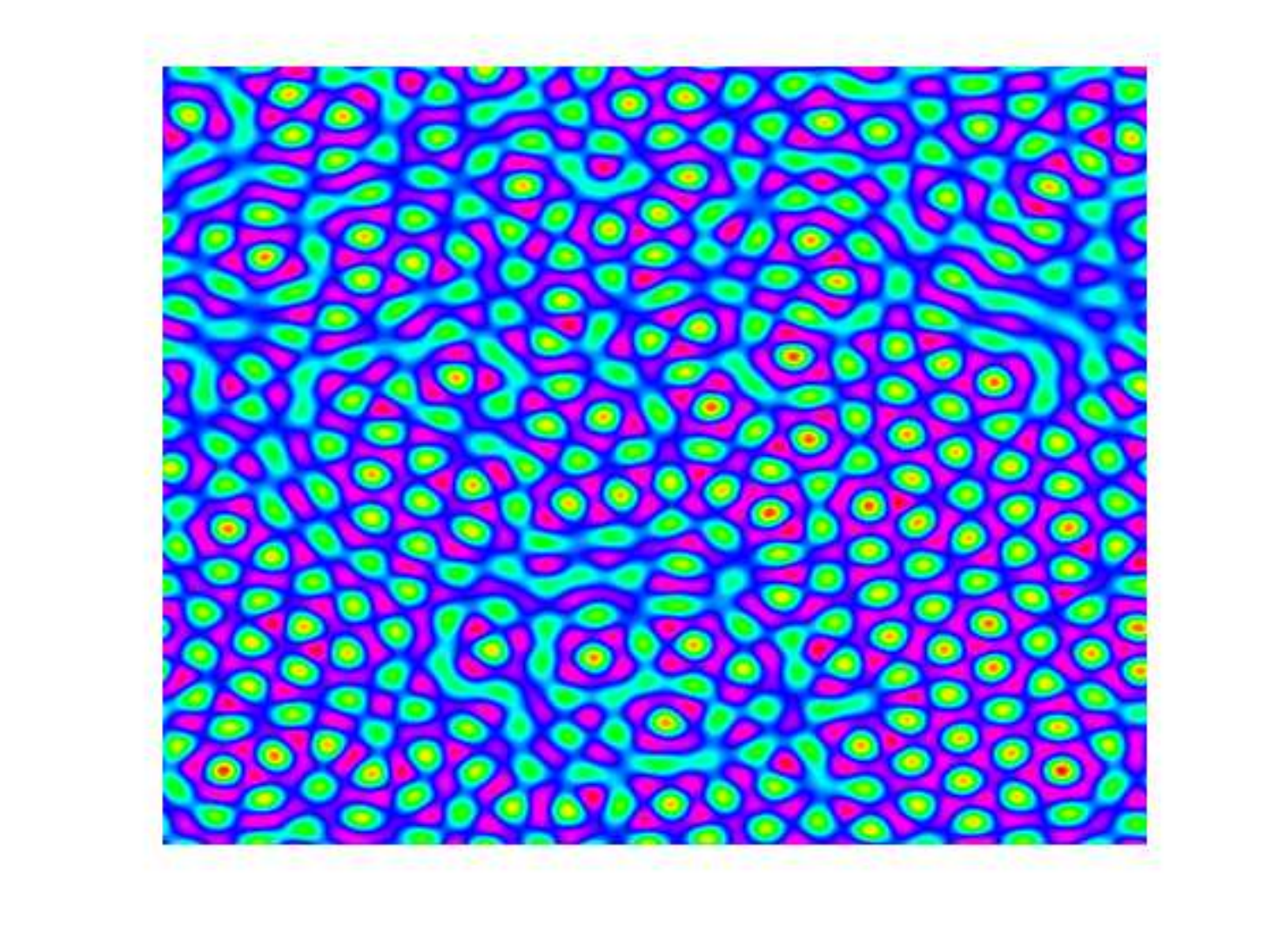}
\includegraphics[width=3.6cm,height=3.6cm]{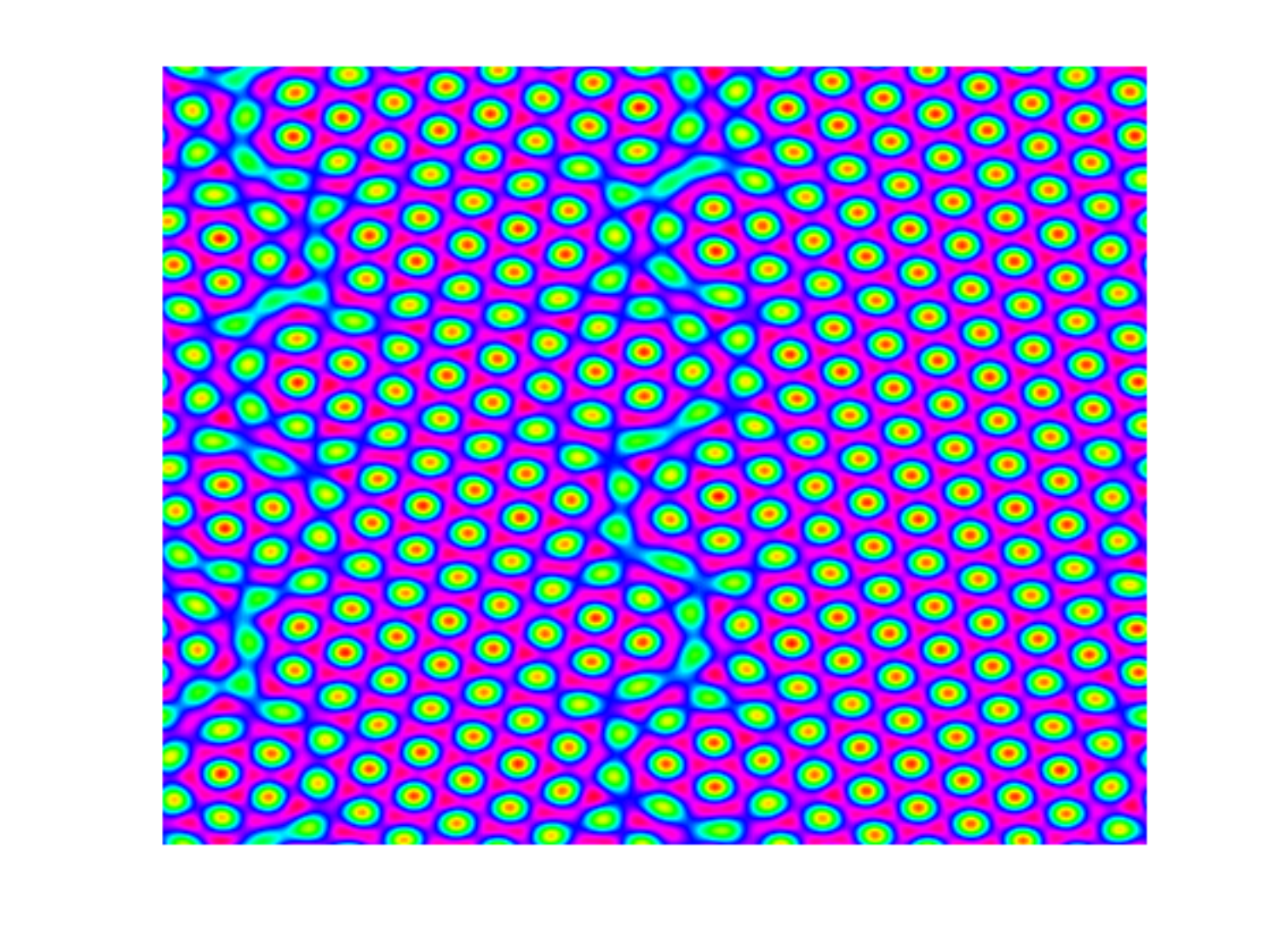}
\includegraphics[width=3.6cm,height=3.6cm]{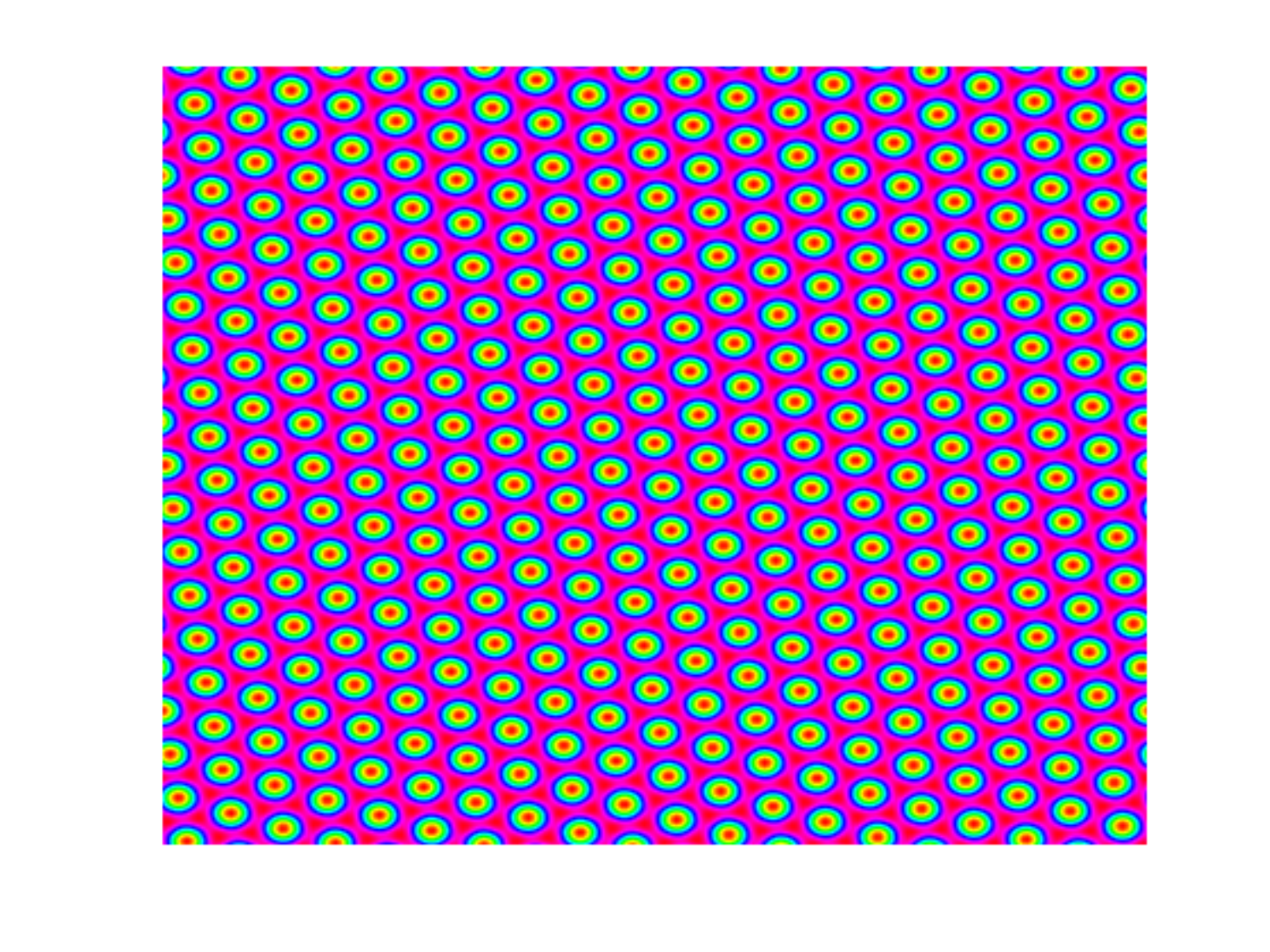}\\
\Xhline{1.2pt}
3S-SAV&\includegraphics[width=3.6cm,height=3.6cm]{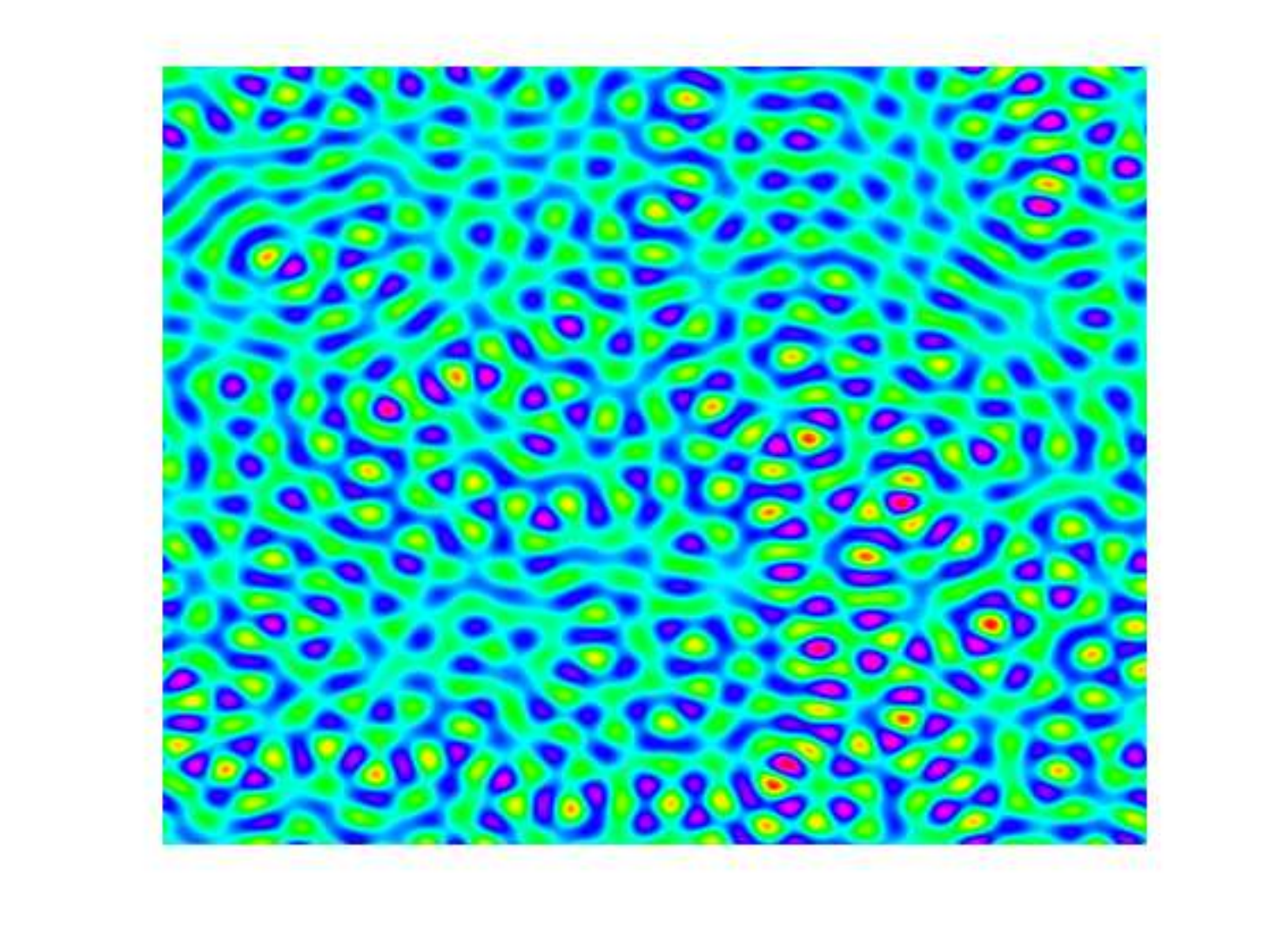}
\includegraphics[width=3.6cm,height=3.6cm]{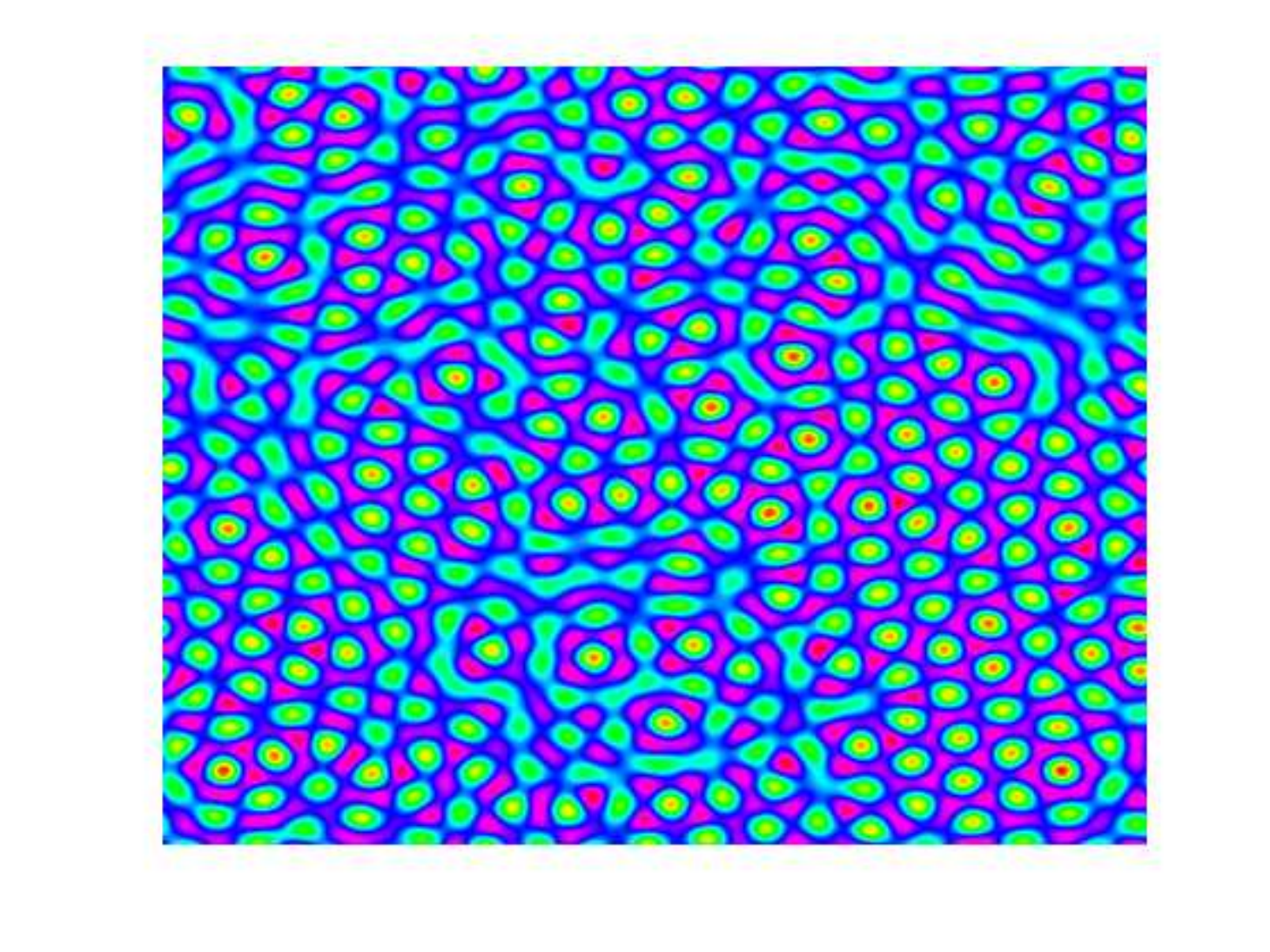}
\includegraphics[width=3.6cm,height=3.6cm]{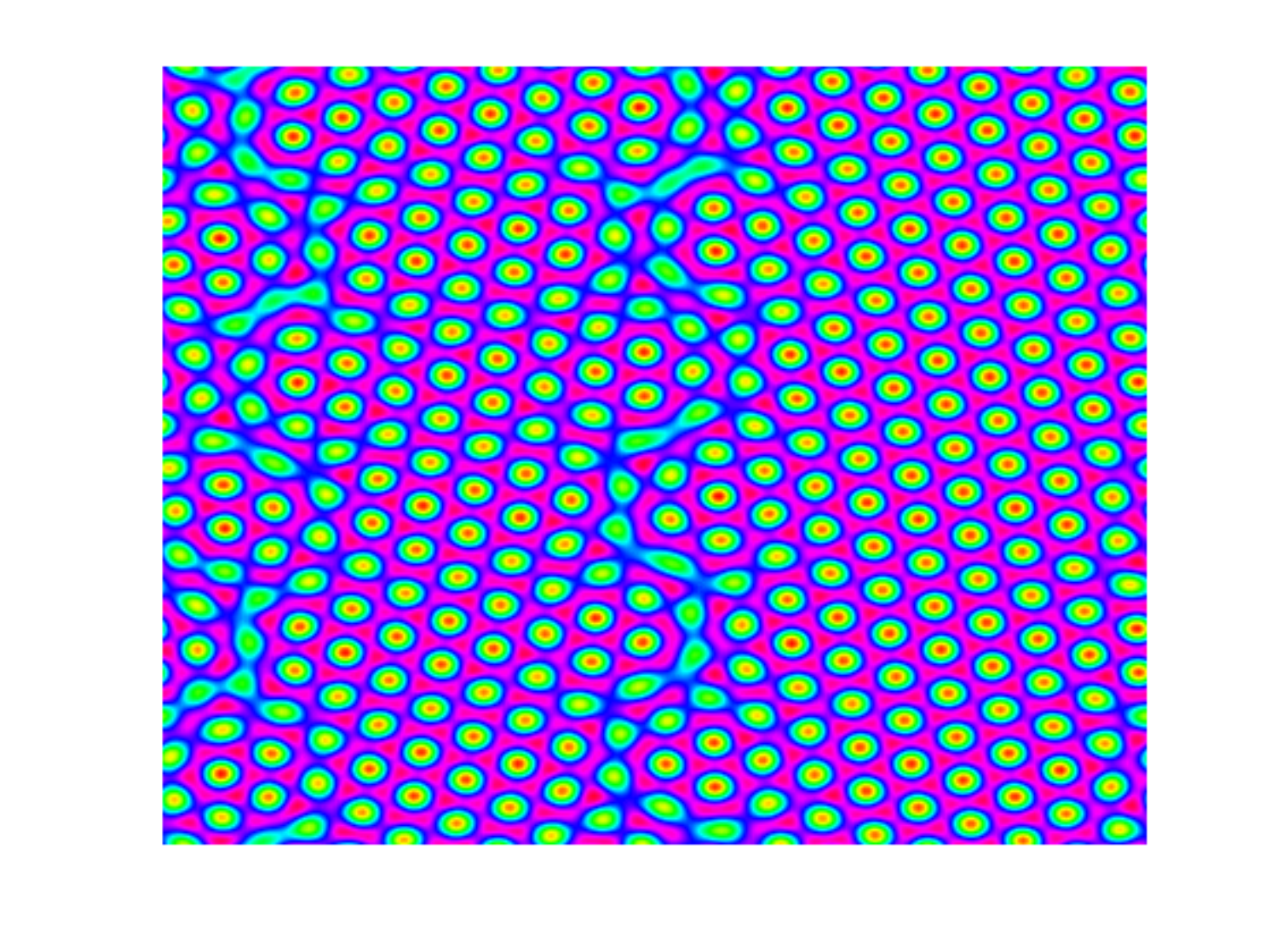}
\includegraphics[width=3.6cm,height=3.6cm]{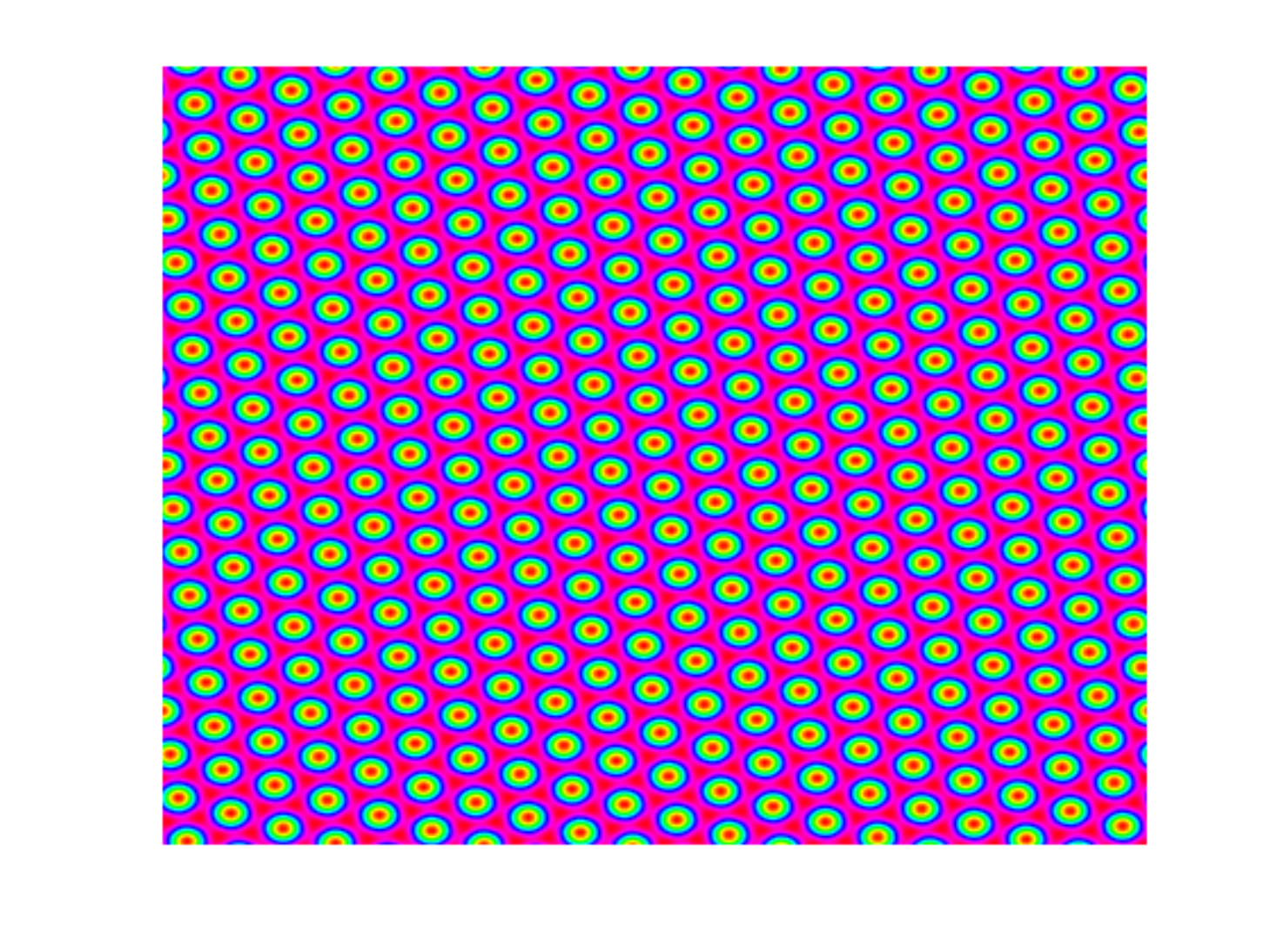}\\
\Xhline{1.2pt}
\end{tabular}
\caption{Configuration evolutions for PFC model by SAV and 3S-SAV schemes are taken at $t=200$, $500$, $1200$, and $6000$.}\label{fig:fig4}
\end{figure}

We show the phase transition behavior of the density field at various times in Figure \ref{fig:fig4}. Similar computation results for phase field crystal model can be found in \cite{ShenA,yang2017linearly}. We investigate the process of crystallization in a supercool liquid by using both SAV and the proposed 3S-SAV schemes. No visible difference is observed.

\section*{Acknowledgement}
No potential conflict of interest was reported by the author. We would like to acknowledge the assistance of volunteers in putting together this example manuscript and supplement.
\bibliographystyle{siamplain}
\bibliography{Reference}

\end{document}